\documentclass[a4j,11pt]{article}

\usepackage[dvipdfmx,hiresbb]{graphicx}
\usepackage[dvipdfmx]{color}
\usepackage{amsmath}
\usepackage{amsthm}
\usepackage{ascmac}
\usepackage{amssymb}
\usepackage{color}
\usepackage{lastpage}
\usepackage{natbib}

\usepackage{geometry}
\numberwithin{equation}{section}
\newtheorem{thm}{Theorem}[section]
\newtheorem{cor}[thm]{Corollary}
\newtheorem{lem}[thm]{Lemma}
\newtheorem{prop}[thm]{Proposition}
\newtheorem{rem}{Remark}

\title{
On Poisson approximations for the Ewens sampling formula when the mutation parameter grows with the sample size}
\author{Koji Tsukuda\footnote{Graduate School of Arts and Sciences, the University of Tokyo, 3-8-1 Komaba, Meguro-ku, Tokyo 153-8902. mail: ctsukuda@g.ecc.u-tokyo.ac.jp}}
\begin{document}
\maketitle

\begin{abstract}
The Ewens sampling formula was firstly introduced in the context of population genetics by Warren John Ewens in 1972, and has appeared in a lot of other scientific fields.
There are abundant approximation results associated with the Ewens sampling formula especially when one of the parameters, the sample size $n$ or the mutation parameter $\theta$ which denotes the scaled mutation rate, tends to infinity while the other is fixed.
By contrast, the case that $\theta$ grows with $n$ has been considered in a relatively small number of works, although this asymptotic setup is also natural.
In this paper, when $\theta$ grows with $n$, we advance the study concerning the asymptotic properties of the total number of alleles and of the counts of components in the allelic partition assuming the Ewens sampling formula, from the viewpoint of Poisson approximations.
\end{abstract}

\newpage
\section{Introduction}
For a positive integer $n$, consider the sequence $\{ C_j^n \}_{j=1}^\infty$ of nonnegative integer-valued random variables satisfying $\sum_{j=1}^n j C_j^n =n$ and $C_j^n =0$ for $j > n$.
For $b=1,\ldots,n$, let us denote $\textbf{C}^n_b=(C_1^n,\ldots,C_b^n)$ and $\textbf{a}_b=(a_1,\ldots,a_b)$.
This $\textbf{C}^n_b$ denotes the component counts in a random combinatorial structure of size $n$.
In the context of population genetics, \cite{RefE} introduced what is called {\it the Ewens sampling formula}
\begin{equation}
{\sf P}(\textbf{C}^n_n = \textbf{a}_n ) = \frac{n!}{(\theta)_n} \prod_{j=1}^n \left( \frac{\theta}{j} \right)^{a_j} \frac{1}{a_j!} 1\left\{ \sum_{j=1}^n j a_j =n \right\}    \label{ESF}
\end{equation}
as the distribution of the allelic partition in a sample of size $n$ from the population which follows the stationary distribution of the infinitely-many neutral allele model with scaled mutation rate $\theta>0$, where $ (\theta)_n$ is the rising factorial $\theta(\theta+1)\cdots(\theta+n-1) $.
See for instance Section 2.5 of \cite{RefF10} for its derivation and basic properties.
Hereafter, we consider \eqref{ESF} as a model of $\{ C^n_j \}_{j=1}^n$.
The unsigned Stirling number of the first kind $\bar{s}(n,k)$ ($k=1,2,\ldots,n$) is the coefficent of $\theta^k$ in $(\theta)_n$, and is in conformity with the number of permutations of $n$ elements with $k$ disjoint cycles.
Hence, if \eqref{ESF} is assumed, the total number $K_n=\sum_{j=1}^n C_j^n$ of alleles included in the sample, in other words the total number of distinct cycles in a random permutation, follows {\it the falling factorial distribution} \citep{RefW1}
\begin{equation}
{\sf P}(K_n=k) = \bar{s}(n,k) \frac{\theta^k}{(\theta)_n}. \label{FFD}
\end{equation}
In this paper, we will present asymptotic properties, especially Poisson approximations, of $\textbf{C}^n_b$ and $K_n$ when both $\theta$ and $n$ increase.

Beyond population genetics domain, the Ewens sampling formula has been widely applied to other fields such as ecology, disclosure risk assessments, nonparametric statistics and so on.
In addition, laws of component counts in a lot of random structures are approximated by the Ewens sampling formula.
For a general review and an up-to-date review with discussions, we refer the reader to Chapter 41 of \cite{RefJKB}, whose write-up was provided by S. Tavar\'e and W.J. Ewens, and to \cite{RefC}, respectively. 
For \eqref{ESF}, \eqref{FFD} and related probabilistic models, a lot of works have discussed asymptotic properties under the situations $n\to\infty$ with fixed $\theta$ or $\theta\to\infty$ with fixed $n$, see for instance \cite{RefF16}.
It is natural to consider some relations between the population size and the sample size.
Since $\theta$ is proportional to the population size in the context of population genetics, \cite{RefF07} and \cite{RefT} discussed the asymptotic behavior of $K_n$ under the settings that both $n$ and $\theta$ simultaneously tend to infinity.
Under this asymptotic setting, \cite{RefF07} established the large deviation principle and \cite{RefT} demonstrated asymptotic properties of the maximum likelihood estimator of $\theta$. 

Following previous works, we set three major goals.
\cite{RefT} extended the asymptotic normality of $K_n$ as $n\to\infty$ with fixed $\theta$, which is due to \cite{RefW}, to the situation when both $\theta$ and $n$ increase.
The first goal of this paper is discussing this result from the viewpoint of Poisson approximations.
Moreover, \cite{RefABT} showed the Poisson process approximation of $\textbf{C}_b^n$ as $n\to\infty$ with fixed $\theta$ when $b$ is fixed or grows with $n$, and \cite{RefAST} established its total variation asymptotics.
The second goal is studying corresponding asymptotic results about $\textbf{C}_b^n$ when $\theta$ grows with $n$.
Furthermore, \cite{RefH} provided a functional central limit theorem for the Ewens sampling formula, and \cite{RefAT} gave its elegant proof via the Poisson process approximation.
Our third goal is to discuss extensions of this and related weak convergence results.

\subsection{Notations}
Consider sequences $\{x_n\}_{n=1}^\infty$ and $\{y_n\}_{n=1}^\infty$.
If $x_n/y_n\to1$, then we write $x_n \sim y_n$.
Let $c <\infty$ be a constant.
If $x_n/y_n \to 0$ then we write $x_n = o(y_n)$, if $x_n/y_n \to c$ then we write $x_n = O(y_n)$, and if $x_n/y_n \to c \neq 0$ then we write  $x_n = \Theta(y_n)$.
Let $\sum_{j=1}^0 x_j =0$ and $\prod_{j=1}^0 x_j =1$ for any sequence $\{ x_\cdot\}$, and let $(x)_0=1$ for any value $x$.
When we consider the limits of $n$ and $\theta$ simultaneously, we use the notation $\lim_{n,\theta}$.

Let $[x^k] f(x)$ denote the coefficient of $x^k$ in the power series expansion of $f(x)$.
Let $f^{(i)}(\cdot)$ denote the $i$-th derivative of  function $f(\cdot)$.
Let $\lfloor \cdot \rfloor$ and $\lceil \cdot \rceil$ denote the floor function and the ceiling function, respectively.
Let $\Gamma(\cdot)$ be the gamma function and $\psi(\cdot) = ( \log \Gamma(\cdot) )'$ the digamma function.
For real $x$, $x^+$ denotes the positive part of $x$.

The space $D[0,1]$ is the set of c\`{a}dl\`{a}g functions on $[0,1]$ endowed with the Skorokhod topology.
The space $L^2 (0,1)$ is equivalence classes of real valued functions on $(0,1)$ which are square integrable with respect to the Lebesgue measure endowed with the $L^2$ topology.

The total variation distance between the laws which random vectors $X$ and $Y$ follow is denoted by $d_{TV}(X,Y)$.
The convergence of $X$ to $Y$ in probability and the weak convergence of $X$ to $Y$ are denoted by $X \to^p Y$ and $X \Rightarrow Y$, respectively.

\subsection{Asymptotic settings}
Letting $c$ be a finite constant, we study the following asymptotic settings in this paper:\\
\qquad Case A: $n/\theta \to \infty$; Case B: $n/\theta \to c>0$; Case C: $n/\theta\to0$;\\
\qquad Case C1: $n/\theta\to0$ and $n^2/\theta\to\infty$; Case C2: $n^2/\theta \to c>0$; Case C3: $n^2/\theta\to0$.\\
This devision is introduced in \cite{RefT}.
It should be noted that in Section 4 of \cite{RefF07}, when $\theta$ does not converge to 0, the relation between $n$ and $\theta$ are divided into Cases A, B, C above and $\theta\to\infty$ with fixed $n$.
Moreover, throughout this paper, we assume that $\theta$ does not decrease as $n$ increase.

\begin{rem}
In Case C3, it holds that $K_n- n \to^p 0$.
Note that when $\theta=o(1/\log{n})$ in which we are not interested since $\theta\to0$, it holds that $K_n-1 \to^p 0$.
These convergences can be checked through showing the convergence in first mean.
\end{rem}

\subsection{Organization}
In Section \ref{sec2}, we review asymptotic results associated with the Ewens sampling formula in the literature which will be discussed in this paper.
Before probabilistic result, in Section \ref{sec3},  let us provide some preliminary evaluations for sequences related to the mean of $K_n$.
Section \ref{sec4} is devoted to show Poisson approximations for $K_n$ and $n-K_n$ in Case A and C, respectively.
Section \ref{sec5} is devoted to discuss independent process approximations for $\textbf{C}_b^n$ in a Ewens partition.
Section \ref{sec6} shows the functional central limit theorems for the Ewens sampling formula when $\theta$ grows with $n$.
In addition, Appendix includes some lemmas used in proofs.

\section{Results in the literature}\label{sec2}
\subsection{Normal and Poisson approximations for $K_n$}\label{sec21}
In the combinatorial context, it is worthwhile to know when typical distributions such as Normal, Poisson or other distributions asymptotically appear.
See for instance \cite{RefFS}.
For the total number $K_n$ of alleles which follows \eqref{FFD}, \cite{RefW} proved the following central limit theorem (CLT for short): For fixed $\theta>0$, 
\begin{equation}
\frac{K_n-\theta\log{n}}{\sqrt{\theta\log{n}}} \Rightarrow N(0,1) \label{WCLT}
\end{equation}
as $n\to\infty$, where $N(0,1)$ is a standard normal variable.
A stronger result, the Poisson approximation for $K_n$, was stated by \cite{RefAT}: For fixed $\theta>0,$
\begin{equation}
d_{TV}(K_n,P_{ {\sf E} [K_n]}) = \Theta \left(\frac{1}{\log{n}}\right)  \label{PAn}
\end{equation}
as $n\to\infty$, where $P_{{\sf E}[K_n]}$ is a Poisson variable with mean ${\sf E}[P_{{\sf E}[K_n]}]={\sf E}[K_n]$.
Later, in order to improve the approximation accuracy, \cite{RefY} provided the following CLT which adopts another standardization: For fixed $\theta>0$, 
\begin{equation}
\frac{K_n-\theta(\log{n}-\psi(\theta))}{\sqrt{\theta(\log{n}-\psi(\theta))}} \Rightarrow N(0,1) \label{YCLT}
\end{equation}
as $n\to\infty$.
Moreover, \cite{RefY} showed the approximation for $K_n$ by a Poisson variable with the approximate mean:  For fixed $\theta>0,$
\begin{equation}
d_{TV}(K_n,P_{\theta (\log{n} - \psi(\theta))}) = O \left(\frac{1}{\log{n}}\right)  \label{YPA}
\end{equation}
as $n\to\infty$, where $P_{\theta (\log{n} - \psi(\theta))}$ is a Poisson variable with mean ${\sf E}[P_{\theta (\log{n} - \psi(\theta))}]=\theta (\log{n} - \psi(\theta))$.

When $\theta$ grows with $n$, the standardization should be changed in many cases.
Let $Z_n = (K_n-\mu)/\sigma$, where $\mu = \theta\log{(1+n/\theta)}$ and $\sigma^2 = \theta (\log{(1+n/\theta)} - n/(n+\theta) )$.
\cite{RefT} showed that 
\begin{equation}
Z_n \Rightarrow 
\begin{cases}
  N(0,1) & {\rm (Case \ A, B, C1) } \\
   (c/2 - P_{c/2})/\sqrt{c/2}  & {\rm (Case \ C2)} \\
   0, & {\rm (Case \ C3)}
\end{cases} \label{TC}
\end{equation}
where $c= \lim_{n,\theta} n^2/\theta$ and $P_{c/2}$ is a Poisson variable with mean ${\sf E}[P_{c/2}]=c/2$.

\begin{rem}
Professor Shuhei Mano pointed out that the proof of Theorem 2 in \cite{RefT} is incorrect in Case C1.
In this remark, let us correct the failure.
As it is stated in the right-hand side in the equation (14) of \cite{RefT}, it holds that 
$ \log {\sf E}[e^{Z_n t}] = -\sigma t + (e^{t/\sigma} -1)\sigma^2 + A + o(1), $
where 
\[ A= \theta \left(\frac{\theta}{n+\theta} -1 \right) \left\{ \frac{t}{\sigma} - (e^{t/\sigma}-1) \right\} - \frac{t}{\sigma}\theta e^{t/\sigma} +(\theta e^{t/\sigma} + n ) \log \left(1+ \frac{\theta}{n+\theta} (e^{t/\sigma}-1) \right) .\]
In Case C1, since $\sigma \sim \sqrt{n^2/2\theta}$, it holds that $\theta \{ n/(n+\theta) \}^3(-t/\sigma)^3 = O(1/\sqrt{\theta})$, and hence
\begin{eqnarray*}
&& (\theta e^{t/\sigma} + n ) \log \left( 1+ \frac{\theta}{n+\theta} (e^{t/\sigma}-1) \right) \\
&=&  (\theta e^{t/\sigma} + n ) \log \left( \frac{n+ \theta e^{t/\sigma}}{n+\theta}  \right)\\
&=& (\theta e^{t/\sigma} + n ) \log \left( e^{t/\sigma} + \frac{n}{n+\theta} (1- e^{t/\sigma}) \right) \\
&=& (\theta e^{t/\sigma} + n ) \left( \frac{t}{\sigma} + \log \left(1+ \frac{n}{n+\theta} (e^{-t/\sigma}-1) \right)  \right) \\
&=&  (\theta e^{t/\sigma} + n ) \left\{ \frac{t}{\sigma} + \frac{n}{n+\theta} (e^{-t/\sigma}-1) -\frac{1}{2}\left(\frac{n}{n+\theta}\right)^2(e^{-t/\sigma}-1)^2 \right\} + o(1).
\end{eqnarray*}
We thus have 
\begin{eqnarray}
A &=& - \frac{\theta n}{n+\theta} \frac{t}{\sigma} + \frac{n}{\sigma} t + \frac{n^2}{n+\theta}(e^{-t/\sigma} -1) - \frac{1}{2} \left( \frac{n}{n+\theta} \right)^2 ( \theta e^{t/\sigma}  +n )  \left( e^{-t/\sigma} -1 \right)^2 + o(1) \nonumber \\
&=& \frac{n}{\sigma} t \left( 1-\frac{\theta}{n+\theta} \right) + n \frac{n}{n+\theta}(e^{-t/\sigma} -1) -\frac{1}{2} \left( \frac{n}{n+\theta} \right)^2 ( \theta e^{t/\sigma}  +n )  (e^{-t/\sigma} -1)^2 + o(1) \nonumber \\
&=& n \left( \frac{n}{n+\theta} \right)  \left( \frac{t}{\sigma} + e^{-t/\sigma} -1 \right) - \frac{n^2 }{2\theta } \frac{ e^{t/\sigma}}{ (1+n/\theta)^2}  \left( 1+ \frac{ne^{-t/\sigma}}{\theta} \right) (e^{-t/\sigma} -1)^2 +o(1) \nonumber \\
&=&  n \left( \frac{n}{n+\theta} \right)  \left( \frac{t}{\sigma} + e^{-t/\sigma} -1 \right)  -  \frac{n^2}{2\theta} (e^{-t/\sigma} -1 )(1-e^{t/\sigma}) \frac{1+n e^{-t/\sigma}/\theta}{(1+n/\theta)^2} + o(1) \label{RRH1}. 
\end{eqnarray}
By using $n/(n+\theta) = n/\theta -n^2/\theta^2 + O(n^3/\theta^3)$, the first term in the right-hand side of \eqref{RRH1} is 
\[ n \left( \frac{n}{\theta } - \frac{n^2}{\theta^2} + O \left( \frac{n^3}{\theta^3} \right) \right)  \left(\frac{t^2}{2\sigma^2} + O \left(\frac{1}{\sigma^3} \right) \right) = \frac{n^2}{\theta } \frac{t^2}{2\sigma^2}   + o(1) \]
The second term in \eqref{RRH1} is also $-n^2t^2/(2\theta\sigma^2) + o(1)$ because it holds that
\[   (e^{-t/\sigma} -1 )(1-e^{t/\sigma}) = \left( -\frac{t}{\sigma} +\frac{t^2}{2\sigma^2} +o \left( \frac{1}{\sigma^2}  \right)\right) \left( -\frac{t}{\sigma} -\frac{t^2}{2\sigma^2} +o \left( \frac{1}{\sigma^2} \right) \right)
= \frac{t^2}{\sigma^2} + O \left( \frac{\theta^2}{n^4} \right) \]
and that $(1+n e^{-t/\sigma}/\theta)/(1+n/\theta)^2 = 1 + O(n/\theta)$.
Therefore $A=o(1)$, and, consequently, $\log {\sf E}[e^{Z_n t}] = - \sigma t +( e^{t/\sigma} -1)\sigma^2  + o(1)\to t^2/2. $
\end{rem}

\begin{rem}
As a corollary to the large deviation principle for $K_n$ when $\theta \to \infty$, \cite{RefF07} provided the following weak law of large numbers in Corollary 4.1: 
\begin{eqnarray}
&& \frac{K_n}{\theta\log{(n/\theta)}} \to^p 1, \quad  {\rm (Case \ A) } \label{FLLN} \\
&& \frac{K_n}{n} \to^p
\begin{cases}
  \log{\left( 1+\frac{1}{c} \right)^c}, & {\rm (Case \ B) } \\
  1,  & {\rm (Case \ C)}  
\end{cases} \nonumber
\end{eqnarray}
and $K_n \to^p n$ as $\theta \to \infty$ with fixed $n$.
These law of large numbers in Cases A, B and C can be obtained directly from the calculation of ${\sf E}[|K_n/{\sf E}[K_n]-1|^2]$, see Proposition 2 of \cite{RefT}.
\end{rem}

\subsection{Independent process approximations for ${\bf C}_b^n$}\label{sec22}
Consider a sequence $\{Z_j\}_{j=1}^\infty$ of independent Poisson variables with ${\sf E}[Z_j] = \theta/j$ for $j=1,2,\ldots$ and denote $\textbf{Z}_b = (Z_1,\ldots,Z_b)$ for a positive integer $b$.
Then, it is well-known that \eqref{ESF} can be derived from {\it the conditioning relation}
\begin{equation}
{\sf P}(\textbf{C}^n_n = \textbf{a}_n ) = {\sf P} \left(  \textbf{Z}_n = \textbf{a}_n \left| \sum_{j=1}^n jZ_j =n \right. \right), \label{CR}
\end{equation}
see for instance \cite{RefW1}.
It means that the dependence in $\{ C_j^n \}_{j=1}^\infty$ is given by the condition $\sum_{j=1}^n j Z_j =n$.
It is of interest to discuss whether the effect of this dependence asymptotically vanishes or not.
It was answered by \cite{RefABT} who showed the small components can be approximated by independent Poisson variables:
For any fixed positive integer $b$, it holds that
\begin{equation}
(C_1^n,\ldots,C_b^n) \Rightarrow (Z_1,\ldots,Z_b)  \label{PA1}
\end{equation}
as $n\to\infty$.
Note that \eqref{PA1} is equivalent to $\lim_{n\to\infty} d_{TV}({\bf C}_b^n, {\bf Z}_b) =0$ because both ${\bf C}_b^n$ and ${\bf Z}_b$ are discrete.

It is more interesting to consider the case that $b$ grows with $n$.
For positive integer $b$, let us denote the total variation distance and the distance in the Wasserstein $\ell^1$ metric between ${\bf C}^n_b=(C_1^n,\ldots,C_b^n)$ and ${\bf Z}_b= (Z_1,\ldots,Z_b)$ by $d_b (n)$ and $d_b^W (n)$, respectively, that is,
\begin{eqnarray*}
d_b (n) &=& d_{TV}( {\bf C}^n_b, {\bf Z}_b ) = \inf_{ {\rm couplings}} {\sf P}( {\bf C}^n_b \neq {\bf Z}_b) ,\\
d_b^W (n) &=& \inf_{ {\rm couplings}} \sum_{j=1}^b {\sf E}[|C_j^n - Z_j|]. 
\end{eqnarray*}
For these quantities, it holds that
\[ d_b(n) = \inf_{ {\rm couplings}} {\sf P}( {\bf C}^n_b \neq {\bf Z}_b) = \inf_{ {\rm couplings}} {\sf P}\left( \sum_{j=1}^b | C_j^n -Z_j | \geq 1 \right) \leq d^W_b(n) . \]
As for the Ewens sampling formula, $d_b^W(n)$ is a convenient measure of approximations because a concrete construction, {\it the Feller coupling}, can be given. 
See \cite{RefABT,RefABT2}.
The Feller coupling is as follows:
Let $\{ \xi_j\}_{j=1}^\infty$ be a sequence of Bernoulli variables with ${\sf P} (\xi_j=1)=p_j=\theta/(\theta+j-1)$ for any $j = 1,2,\ldots$.
Then, the Ewens sampling formula \eqref{ESF} is given as the joint distribution of
\[
C_1^n = \sum_{i=1}^{n-1} \xi_i\xi_{i+1} + \xi_n
\]
and
\[
C_j^n = \sum_{i=1}^{n-j} \xi_i (1-\xi_{i+1})\cdots (1-\xi_{i+j-1}) \xi_{i+j} + \xi_{n-j+1}(1-\xi_{n-j+2})\cdots(1-\xi_{n}) 
\]
for $j=2,\ldots,n$.
Moreover, define
\begin{eqnarray*}
C_j^\infty = \sum_{i=1}^{\infty} \xi_i (1-\xi_{i+1})\cdots (1-\xi_{i+j-1}) \xi_{i+j}
\end{eqnarray*}
for $j=1,\ldots,n$, then $C_j^\infty$ follows the independent Poisson distribution with mean ${\sf E}[C_j^\infty] =\theta/j$ for any $j=1,2,\ldots$.
That is because the convergences in probability $\xi_n \to^p 0$ and
$ \xi_{n-j+1}(1-\xi_{n-j+2})\cdots(1-\xi_n) \to^p 0 $
for any $j=2,3,\ldots$ yield that $C_j^n \Rightarrow C_j^\infty$, and so \eqref{PA1} yields that $C_j^\infty =^d Z_j$ for any $j=1,2,\ldots$.
By using this construction, \cite{RefABT} proved the Poisson process approximation for $b$ growing with $n$: 
\begin{eqnarray}
&& d_b(n)  \to 0 \Leftrightarrow b=o(n); \\
&& d_b(n)  \leq \frac{b\theta}{\theta+n} \left( \theta + \frac{n}{\theta+n-b} \right); \label{tbA} \\
&& d_b^W(n) \leq \frac{b\theta}{\theta+n-b}\left(\theta+\frac{n}{\theta+n}\right); \label{wbA} \\
&& d_n^W(n) = O(1); \label{ub} 
\end{eqnarray}
if $\theta\geq 1$ then
\begin{equation}
 \frac{\theta(\theta-1)b}{\theta+n-1} \left\{ 1 - \frac{(\theta-1) (b+1)}{4(\theta+n-1)} \right\} \leq d_b^W(n) \leq \frac{b\theta(\theta+1)}{\theta+n}. \label{wb1}
\end{equation}
Note that \eqref{tbA}, \eqref{wbA} and \eqref{wb1} are not asymptotic results.
Lower bound results for the total variation distance, which complement \eqref{tbA}, were given by \cite{RefABT}: $\liminf_{n\to\infty} n d_b(n) \geq (b\theta|\theta-1|/2) \exp\left(-\theta\sum_{j=1}^b 1/j \right)$; and by \cite{RefB}: if $\theta\neq 1$ then $d_b(n) \geq c_3 b/n$ for some $c_3=c_3(\theta)>0$.

Another compelling result for evaluating $d_b(n)$ is deriving the leading term of $d_b(n)$, which were given by \cite{RefAST} for general logarithmic assemblies.
If the Ewens sampling formula is considered, the statement is as follows:
If $b=o(n/\log{n})$ then
\begin{equation}
d_b(n) = \frac{|1-\theta|}{2n}{\sf E} [|T_{0b} - \theta b|] + o\left(\frac{b}{n}\right),  \label{AST}
\end{equation}
where $T_{0b}= \sum_{j=1}^b j Z_j$.
As it is stated in Corollary 4 of their paper, if $\theta\neq 1$ and if $b=o(n/\log{n})$ then the leading term of  $d_b(n)$ is given by the first term in the right-hand side of \eqref{AST}.

\subsection{Functional central limit theorems}\label{fcs}
The results by \cite{RefABT} provide an elegant way to derive asymptotic properties.
Among others, by using \eqref{ub}, \cite{RefAT} provided an alternative proof of the functional central limit theorem for the Ewens sampling formula which was originally proven by \cite{RefH}:
The random process
\begin{equation}
 X^1_n(\cdot) = \left( \frac{ \sum_{i=1}^{\lfloor n^u \rfloor} C_j^n - u \theta \log{n} }{ \sqrt{\theta \log{n}} } \right)_{0\leq u\leq1}  \label{RF0}
\end{equation}
converges weakly to $(B(u))_{0 \leq u\leq1}$ in $D[0,1]$ as $n\to\infty$, where $B(\cdot)$ is a standard Brownian motion.
This approach is generalized to broader logarithmic structures.
See \cite{RefAST} and \cite{RefABT1}.
Moreover, by using the Poisson process approximation, \cite{RefT1} provided a weighted version in $L^2(0,1)$:
Both of the random processes
\begin{equation}
 X^2_n(\cdot) = \left( \frac{ \sum_{i=1}^{\lfloor n^u \rfloor} C_j^n -  \theta \sum_{j=1}^{\lfloor n^u \rfloor} 1/j  }{ \sqrt{ \theta \sum_{j=1}^{\lfloor n^u \rfloor} 1/j } }  \right)_{0 < u < 1} \label{RF1}
\end{equation}
and
\begin{equation}
 X^3_n(\cdot) = \left( \frac{ \sum_{i=1}^{\lfloor n^u \rfloor} C_j^n - u \theta \log{n} }{ \sqrt{u \theta \log{n}} } 1 \left\{u > \frac{\varepsilon}{\log{n} } \right\}  \right)_{0 < u < 1} \label{RF2}
\end{equation}
converge weakly to $(B(u)/\sqrt{u})_{0<u<1}$ in $L^2(0,1)$ as $n\to\infty$, where $\varepsilon$ is a positive constant.

\begin{rem}
In the case that $\theta=1$, the weak convergence of $X_n^1(\cdot)$ in $D[0,1]$ was provided by \cite{RefDP}.
\end{rem}

Let $R_j$ be the $j$-th cycle length in a random permutation of $n$ which has $K_n$ disjoint cycles, and the {\it loglength} of $j$-th cycle is defined by $\log_n{R_j}$.
Consider its empirical distribution function
\[ F_n (u) = \frac{ \sum_{j=1}^{K_n} 1\{ \log_n{R_j} \leq u \} }{K_n} = \frac{ \sum_{j=1}^{\lfloor n^u \rfloor} C_j^n }{K_n} \]
for $0 \leq u \leq 1$.
Define the random processes
\begin{eqnarray}
X^4_n(\cdot)  &=& \left ( \sqrt{\theta \log{n}} (F_n (u) - u) \right)_{0\leq u \leq 1},  \label{EB} \\
X^5_n(\cdot) &=& \left ( \sqrt{\theta \log{n}} \frac{(F_n (u) - u)}{\sqrt{u(1-u)}} 1\left\{ \frac{\varepsilon}{\log{n}} < u < 1- \frac{\varepsilon}{\log{n}} \right\} \right)_{0 < u < 1}  \label{EBT} ,
\end{eqnarray}
where $\varepsilon$ is a positive constant.
When $\theta=1$, the weak convergence of $X_n^4 (\cdot)$ to a standard Brownian bridge  $(B^\circ (u))_{0\leq u \leq 1}$ in $D[0,1]$ was shown by \cite{RefDP}, see Notes (2) after Theorem in their paper.
Its extension to the Ewens sampling formula and $L^2$ version are presented as follows, which may have not appeared in the literature.

\begin{prop}\label{LR0}
(i) The random process $X^4_n(\cdot)$ converges weakly to  $(B^\circ (u))_{0\leq u \leq 1}$ in $D[0,1]$ as $n\to\infty$. \\
(ii) The random process $X^5_n(\cdot)$ converges weakly to $(B^\circ (u)/\sqrt{u(1-u)})_{0\leq u \leq 1}$ in $L^2(0,1)$ as $n\to\infty$.
\end{prop}

We omit its proof  because we will present an extended version in Proposition \ref{LRC}.
From Proposition \ref{LR0}, it follows from the continuous mapping theorem that
\[
\sup_{0\leq u \leq 1} |X^4_n(u)| \Rightarrow \sup_{0\leq u \leq 1} |B^\circ(u)|,  \quad
\int_{0}^{1} |X^5_n(u)|^2 du \Rightarrow \int_{0}^{1}  \frac{ \left| B^\circ(u) \right|^2}{u(1-u)} du
\]
as $n\to\infty$.

\begin{rem}
As it is stated in \cite{RefDP}, Proposition \ref{LR0} means that $F_n(u)$ is nearly $u$, which is the distribution function of the standard uniform distribution.
\end{rem}

\subsection{Auxiliary results}\label{AUX}
In this subsection, let us set out some auxiliary results concerning Poisson approximations which will be used in the proofs of our statements.

Consider a sequence of independent Bernoulli variables $\{ \xi_j\}_{j=1}^\infty$ and its partial sum $S_n=\sum_{j=1}^n \xi_j$, where ${\sf P}(\xi_j=1)=p_j$ for any $j=1,2,\ldots$.
Then, by using the Chen--Stein method, Theorems 1 and 2 of \cite{RefBH} gave the sharp bound for the Poisson approximation for a partial sum of Bernoulli variables:
For a Poisson variable $P_\lambda$ with mean ${\sf E}[P_\lambda] = \lambda=\sum_{j=1}^n p_j$, it holds that
\begin{equation}
\frac{1\wedge\lambda^{-1}}{32} \sum_{j=1}^n p_j^2 \leq d_{TV}(S_n , P_\lambda) \leq \frac{1-e^{-\lambda}}{\lambda} \sum_{j=1}^n p_j^2.   \label{BH}
\end{equation}
Moreover, from a property of the Hellinger integral, a bound for the total variation distance between two Poisson distributions were given in Theorem 2.1 of \cite{RefYan}: For Poisson variables $P^1_{\lambda_1}$ and $P^2_{\lambda_2}$ with respective means $\lambda_1$ and $\lambda_2$, it holds that
\begin{equation}
 d_{TV}(P^1_{\lambda_1}, P^2_{\lambda_2})  \leq \min \left( |\sqrt{\lambda_1}-\sqrt{\lambda_2}|,|\lambda_1-\lambda_2| \right) . \label{Yan}
\end{equation}

%%%%%%%%%%%%%%%%%%%%%%%%%%%%%%%%%%%%
\section{Preliminary results}\label{sec3}
Before discussing probabilistic results, let us show asymptotic evaluations on sums of sequences which will be used.
Consider two sequences $\{p_j\}_{j\geq1}$ and $\{q_j\}_{j\geq1}$ given by
\begin{equation}
p_j = \frac{\theta}{\theta+j-1}, \quad q_j=1-p_j= \frac{j-1}{\theta+j-1}. \quad (j=1,2,\ldots,n)  \label{pq}
\end{equation}

\begin{prop}\label{Prop1}
(i) It holds that
\begin{equation}
\frac{n}{2(\theta+n)} \leq \sum_{j=1}^n p_j - \theta \log{\left(1+ \frac{n}{\theta} \right)} \leq \frac{n}{\theta+n}  \label{P1i}
\end{equation}
and that
\begin{equation}
 0 \leq \sum_{j=1}^n p_j^2 - \frac{n\theta}{n+\theta} \leq 1. \label{P1ii}
\end{equation}
Especially, in Case A, it holds that
\[ \sum_{j=1}^n p_j \sim \theta \log\left(1+\frac{n}{\theta}\right), \]
and that if $\theta\to\infty$ then
\[ \sum_{j=1}^n p_j^2 \sim \theta. \]
(ii) It holds that
\begin{equation} 
 \frac{n(n-1)}{2(\theta+n)} \leq \sum_{j=1}^n q_j \leq \frac{n(n-1)}{2\theta}  \label{P1iii}
\end{equation}
and that
\begin{equation} 
  \frac{n(n-1)(2n-1)}{6(\theta+n)^2} \leq \sum_{j=1}^n q_j^2 \leq \frac{n(n-1)(2n-1)}{6\theta^2}.   \label{P1iv}
\end{equation}
Especially, in Case C, it holds that
\[ \sum_{j=1}^n q_j \sim \frac{n^2}{2\theta}, \quad \sum_{j=1}^n q_j^2 \sim \frac{n^3}{3\theta^2}. \]
\end{prop}

Proof.
(i) Since
\[ \frac{1}{2} \sum_{j=1}^n\left(  \frac{1}{\theta+j-1} -\frac{1}{\theta+j}\right) \leq  \sum_{j=1}^n\frac{1}{\theta+j-1} - \int_{\theta}^{\theta+n} \frac{dx}{x}  \leq \sum_{j=1}^n \left(\frac{1}{\theta+j-1} -\frac{1}{\theta+j}\right),  \]
the result \eqref{P1i} holds.
Since 
\[ \int_{\theta}^{\theta+n} \frac{dx}{x^2} \leq \sum_{j=1}^n \frac{1}{(\theta+j-1)^2} \leq \frac{1}{\theta^2} + \int_{\theta}^{\theta+n} \frac{dx}{x^2}, \]
the result \eqref{P1ii} holds. 

(ii) From
\[ \frac{j-1}{\theta+n} \leq q_j \leq \frac{j-1}{\theta}\]
and from
\[ \frac{(j-1)^2}{(\theta+n)^2} \leq q_j^2 \leq \frac{(j-1)^2}{\theta^2}, \]
the results \eqref{P1iii} and \eqref{P1iv} follow.
This completes the proof.

\begin{prop}\label{P2}
In Case A, it holds that
\begin{equation} 
\sum_{j=1}^n p_j - \theta(\log{n} - \psi(\theta)) = O\left( \frac{\theta^2}{n} \right). \label{YB}
\end{equation}
\end{prop}

Proof.
It follows from 
\begin{equation}
\sum_{j=1}^n p_j  = \theta (\psi(n+\theta)  - \psi(\theta))  \label{meanr}
\end{equation}
that the left-hand side of \eqref{YB} is $\theta (\psi(n+\theta)  - \log{n})$.
Since 
\[\psi(n+\theta) = -\gamma - \frac{1}{n+\theta} + (n+\theta) \sum_{j=1}^\infty \frac{1}{j(n+\theta+j)} \]
and $ \log{n} =\sum_{j=1}^n 1/j - \gamma + O(1/n)$ as $n\to\infty$, it holds that
\begin{eqnarray*}
\psi(n+\theta)- \log{n} &=& -\sum_{j=1}^n \frac{1}{j}  + \sum_{j=1}^\infty\left( \frac{1}{j} - \frac{1}{n+\theta+j} \right) + O\left(\frac{1}{n}\right) \\
&\leq & -\sum_{j=1}^n \frac{1}{j}  + \sum_{j=1}^\infty\left( \frac{1}{j} - \frac{1}{n+ \lfloor \theta \rfloor +1 +j} \right) + O\left(\frac{1}{n}\right) \\
&=& \sum_{j=n+1}^{n+\lfloor \theta \rfloor }  \frac{1}{j}  + O\left(\frac{1}{n}\right) .
\end{eqnarray*}
In Case A, the first term in the right-hand side is $O(\theta/n)$.
This completes the proof.

\begin{rem}
As it is stated in \eqref{YCLT}, \cite{RefY} discussed the asymptotic normality of $K_n$ standardized by $\theta(\log{n}- \psi(\theta))$, which means that $\psi(n+\theta)$ is approximated by $\log{n}$ from \eqref{meanr}.
If $\theta^{3}/(n^2 \log(n/\theta)) \to \infty$, the bound in \eqref{YB} is meaningless to discuss CLT.
On the other hand, if $\theta^2/n\to0$ the centering by $\theta(\log{n}-\psi(\theta))$ is better than centering by $\theta\log(1+n/\theta)$, which was used in Corollary 2 of \cite{RefT}, because $\sum_{j=1}^n p_j - \theta\log(1+n/\theta) = \Theta\left( n/(n+\theta) \right)$.
\end{rem}

\begin{prop}\label{P3B}
In Case C, it holds that
\[ \sum_{j=1}^n p_j  - \sum_{j=1}^n \frac{\theta}{j}\left( \frac{n}{\theta} \right)^j  =O\left( \frac{n^2}{\theta} \right) . \]
\end{prop}

Proof.
The triangle inequality yields that
\begin{eqnarray*}
&& \left| \sum_{j=1}^n \left\{ p_j  - \frac{\theta}{j}\left( \frac{n}{\theta} \right)^j \right\}  \right| \\
&\leq&  \left| \sum_{j=1}^n p_j  - \theta\log\left(1+\frac{n}{\theta} \right) \right| + \left| \theta\log\left(1+\frac{n}{\theta} \right) - \log\left( 1 -\frac{n}{\theta} \right)^{-\theta}   \right| \\ 
&&   + \left| \log\left( 1 -\frac{n}{\theta} \right)^{-\theta} - \sum_{j=1}^n  \frac{\theta}{j}\left( \frac{n}{\theta} \right)^j  \right|  .
\end{eqnarray*}
The first term is $O(n/(n+\theta))=O(n/\theta)$, the second term is $\theta \log (1- n^2/\theta^2) = O(n^2/\theta)$, and from $\log{(1-x)}^{-1}=x+x^2/2+\cdots$ as $x\to 0$ the third term is 
\[ \sum_{j=n+1}^\infty \frac{\theta}{j} \left( \frac{n}{\theta} \right)^j
\leq \frac{\theta}{n}\sum_{j=n+1}^\infty \left(\frac{n}{\theta}\right)^j
=\frac{\theta}{n}\frac{(n/\theta)^{n+1}}{1-n/\theta} 
=O\left(\left( \frac{n}{\theta} \right)^n \right). \]
This completes the proof.

%%%%%%%%%%%%%%%%%%%%%%%%%%%%%%%%%%%%
\section{Poisson approximations for the total number of alleles}\label{sec4}

Introduce two Poisson variables $Z^{A}$ and $Z^{C}$ whose means are given by $\lambda_A= {\sf E}[K_n] = \sum_{j=1}^n p_j$ and $\lambda_C = n- {\sf E}[K_n] = \sum_{j=1}^n q_j$, respectively, where $K_n$ follows \eqref{FFD} and where $\{p_j\}_{j=1}^\infty$ and $\{q_j\}_{j=1}^\infty$ are given in \eqref{pq}.
Poisson approximations corresponding to \eqref{PAn} are given in the following proposition.

\begin{prop}\label{PAP}
(i) In Case A,
\begin{equation}
d_{TV}(K_n, Z^{A}) \leq \frac{n\theta + n+\theta}{\theta(n+\theta) \log(1+ n/\theta) + n/2 },  \label{PAP1}
\end{equation}
and 
\[ d_{TV}(K_n, Z^{A}) = \Theta\left( \frac{1}{\log(n/\theta)} \right). \]
(ii) In Case C, 
\begin{equation}
d_{TV}(n-K_n, Z^{C}) \leq \frac{2n(n+\theta)}{3\theta^2 } (1-e^{-n^2/2\theta}), \label{PAP2}
\end{equation}
and 
\[ d_{TV}(n-K_n, Z^{C}) = \Theta\left( \frac{n}{\theta} \left(1 \wedge \frac{n^2}{\theta} \right) \right). \]
\end{prop}

Proof.
Let $\{\xi_j \}_{j =1}^\infty$ and $\{\zeta_j \}_{j=1}^\infty$ be sequences of Bernoulli variables with respective parameters ${\sf P}(\xi_j=1) = p_j $ and  ${\sf P}(\zeta_j=1) = q_j $ for $j=1,2,\ldots$.
Then, it holds that $K_n =^d \sum_{i=1}^n \xi_i $ and that $n- K_n =^d \sum_{i=1}^n \zeta_i $.
To prove the desired results, we will use \eqref{BH} and Proposition \ref{Prop1}.

(i) The result \eqref{PAP1} follows from
\[ d_{TV}(K_n, Z^{A}) \leq  \frac{\sum_{j=1}^n p_j^2}{\sum_{j=1}^n p_j} \leq \frac{n\theta/(n+\theta) +1}{\theta\log(1+n/\theta) + n/\{ 2(\theta+n) \} }. \]
Since $\sum_{j=1}^n p_j \to \infty$, it holds that
\[  d_{TV}(K_n, Z^{A}) \geq  \frac{\sum_{j=1}^n p_j^2}{32\sum_{j=1}^n p_j} \geq \frac{n\theta/(n+\theta) }{32 \theta\log(1+n/\theta) } \]
for enough large $n$.
Two displays above imply $ d_{TV}(K_n, Z^{A}) = \Theta\left( 1/ \log(n/\theta) \right)$.

(ii) The result \eqref{PAP2} follows from
\begin{eqnarray*}
d_{TV}(n-K_n, Z^{C}) &\leq& (1-e^{-\sum_{j=1}^n q_j}) \frac{\sum_{j=1}^n q_j^2}{\sum_{j=1}^n q_j} \\
&\leq&  (1-e^{-n^2/2\theta}) \frac{n(n-1)(2n-1)/(6\theta^2)}{n(n-1)/ \{ 2(\theta+n) \} }. 
\end{eqnarray*}
In Case C1, since $\sum_{j=1}^n q_j \to \infty$, it holds that
\[  d_{TV}(n- K_n, Z^{C}) \geq  \frac{\sum_{j=1}^n q_j^2}{32\sum_{j=1}^n q_j} \geq \frac{n(n-1)(2n-1)/ \{6(\theta+n)^2 \} }{ 32 n(n-1)/ 2\theta } \]
for enough large $n$.
Two displays above imply the $ d_{TV}(n- K_n, Z^{C}) = \Theta\left( n/ \theta \right)$.
In Case C2, since $1-e^{-\sum_{j=1}^n q_j} \leq 1$ and since $(1/\sum_{j=1}^n q_j) $ is bounded by some constant for enough large $n$, the same evaluation  provides $ d_{TV}(n- K_n, Z^{C}) = \Theta\left( n/ \theta \right)$.
In Case C3, since $\sum_{j=1}^n q_j \to 0$, it holds that $1-e^{-\sum_{j=1}^n q_j} \sim n^2/(2\theta)$ and that 
\[ d_{TV}(n-K_n,Z^C) \geq \frac{1}{32} \sum_{j=1}^n q_j^2 \geq \frac{n(n-1)(2n-1)}{192 (\theta+n)^2} \]
for enough large $n$.
We thus have $ d_{TV}(n-K_n, Z^{C}) = \Theta\left( n^3/ \theta^2 \right)$.
This completes the proof.

\begin{rem}
From asymptotic properties of the Poisson distribution and Proposition \ref{PAP}, the result of \eqref{TC} in Cases A and C can be derived.
\end{rem}

In Proposition \ref{PAP}, we have considered Poisson variables with rigorous means ${\sf E}[K_n]$ and $n - {\sf E}[K_n]$.
Next, let us discuss centerings by approximate means presented by \cite{RefY} and \cite{RefT} from the viewpoint of Poisson approximation.
Introduce three Poisson variables $Y^{A}, Y^a$ and $Y^C$ whose means are given by $\mu_A= {\sf E}[Y^{A}]= \theta \log(1+n/\theta)$, $\mu_a={\sf E} [Y^a] = \theta (\log{n}- \psi(\theta))$ and $\mu_C = {\sf E}[Y^C] = n- \theta \log(1+n/\theta),$ respectively.

\begin{lem}
(i) In Case A, it holds that
\begin{equation}
d_{TV}(Z^{A}, Y^A) = O\left( \frac{1}{\sqrt{\theta\log(n/\theta)}} \right)  \label{L51}
\end{equation}
and that
\begin{equation}
 d_{TV}(Z^{A}, Y^a) = O\left( \frac{\theta^{3/2}}{n \sqrt{\log(n/\theta)}} \right).   \label{L52}
\end{equation}
(ii) In Case C, it holds that
\begin{equation}
d_{TV}(Z^{C}, Y^C) = O\left( \frac{1}{\sqrt{\theta}} \left( 1 \wedge \frac{n}{\sqrt{\theta}} \right) \right).  \label{L53}
\end{equation}
\end{lem}

Proof.
We will use \eqref{Yan}.
(i) First we see \eqref{L51}.
Since $\lambda_A$ and $\mu_A$ tend to infinity in Case A, $  |\sqrt{\lambda_A} - \sqrt{\mu_A} |  \leq |\lambda_A -\mu_A |$ for enough large $n$.
Moreover, by using Proposition \ref{Prop1}, it holds that
\[ |\sqrt{\lambda_A} - \sqrt{\mu_A} |   = \frac{ |\lambda_A - \mu_A| }{\sqrt{\lambda_A} + \sqrt{\mu_A} } =   \frac{ |\sum_{j=1}^n p_j  - \theta \log(1+n/\theta) | }{ \sqrt{\sum_{j=1}^n p_j} + \sqrt{\theta \log(1+n/\theta)}  } = O\left( \frac{n /  (n + \theta) }{ \sqrt{\theta \log(1+n/\theta)}} \right), \]
and hence \eqref{L51}.

Next we see \eqref{L52}.
By using Propositions \ref{Prop1} and \ref{P2}, it holds that
\[|\sqrt{\lambda_A} - \sqrt{\mu_A} |  =  \frac{ |\lambda_A - \mu_a| }{\sqrt{\lambda_A} + \sqrt{\mu_a} } =   \frac{ | \sum_{j=1}^n p_j - \theta (\log{n} - \psi(\theta)) |}{\sqrt{\sum_{j=1}^n p_j} + \sqrt{\theta (\log{n} - \psi(\theta))}  } = O\left( \frac{\theta^{2} /n }{ \sqrt{\theta \log(1+n/\theta)}} \right), \]
and hence \eqref{L52}.

(ii) First consider Case C1.
Since $\lambda_C$ and $\mu_C$ tend to infinity in Case C1, $|\sqrt{\lambda_C} - \sqrt{\mu_C} |  \leq |\lambda_C -\mu_C |$ for enough large $n$. 
By using Proposition \ref{Prop1}, it holds that
\[ |\sqrt{\lambda_C} - \sqrt{\mu_C} |  = \frac{ |\lambda_C - \mu_C| }{\sqrt{\lambda_C} + \sqrt{\mu_C} } =   \frac{ |\sum_{j=1}^n q_j  - (n - \theta \log(1+n/\theta)) | }{ \sqrt{\sum_{j=1}^n q_j} + \sqrt{n - \theta \log(1+n/\theta)}  } = O\left( \frac{n/(\theta + n) }{ \sqrt{ n^2/\theta }} \right) = O\left( \frac{1}{\sqrt{\theta}} \right), \]
and hence \eqref{L53} holds as $d_{TV}(Z^C,Y^C)=O(1/\sqrt{\theta})$.

Next consider Case C2.
The magnitude relationship of $|\sqrt{\lambda_C} - \sqrt{\mu_C} |$ and $|\lambda_C -\mu_C |$ is not determined, but they have the same bound $O(1/\sqrt{\theta})$ because $1/\sqrt{\theta} = \Theta(1/n)$.
Hence \eqref{L53} holds as $d_{TV}(Z^C,Y^C)=O(1/\sqrt{\theta})$.

Finally, consider Case C3.
Since $\lambda_C$ and $\mu_C$ tend to 0, $|\sqrt{\lambda_C} - \sqrt{\mu_C} |  \geq |\lambda_C -\mu_C |$ for enough large $n$. 
By using Proposition \ref{Prop1}, it holds that
\[  |\lambda_C - \mu_C|  =   \left|  \sum_{j=1}^n q_j  - \left( n- \theta \log\left( 1+ \frac{n}{\theta} \right) \right) \right|  = O\left( \frac{n}{ \theta} \right), \]
and hence \eqref{L53} holds as $d_{TV}(Z^C,Y^C)=O(n/\theta)$.
This completes the proof.
\\

From what has already been proven, the triangle inequality yields the following Poisson approximations corresponding to \eqref{YPA}.

\begin{prop}
(i) In Case A, if $(\log{(n/\theta)})/\theta \to \infty$ then
\[ d_{TV}(K_n, Y^A) = O\left(\frac{1}{\sqrt{\theta\log(n/\theta)}} \right), \]
and if $(\log{(n/\theta)})/\theta = O(1)$ then
\[ d_{TV}(K_n, Y^A) = O\left(\frac{1}{\log(n/\theta)} \right). \]
Moreover, in Case A, if $(\theta^{3}\log{(n/\theta)})/n^2 =O(1)$ then
\[ d_{TV}(K_n, Y^a) = O\left(\frac{1}{\log(n/\theta)} \right), \]
and if $(\theta^{3}\log{(n/\theta)})/n^2 \to \infty$ and $\theta^{3}/(n^2\log{(n/\theta)}) =O(1)$ then
\[
d_{TV}(K_n, Y^a) = O\left( \frac{\theta^{3/2}}{n \sqrt{\log(n/\theta)}} \right). 
\]
(ii) In Case C, it holds that
\[ d_{TV}(n- K_n, Y^C) = O\left( \frac{n}{\theta} \right). \]
\end{prop}

%%%%%%%%%%%%%%%%%%%%%%%%%%%%%%%%%%%%
\section{On independent process approximation of component counts}\label{sec5}

\subsection{Case A}

First we see the asymptotic independence of small components ${\bf C}^n_b = (C_1^n,\ldots, C_b^n)$ in Case A when $\theta \geq 1$ for some $n$, recalling that we assume that $\theta$ does not decrease as $n$ increase.
We will not discuss the other case, $\theta<1$ for all $n$, because we are interested in large $\theta$

Consider $\{ Z_j \}_{j=1}^\infty$ defined in Subsection 2.2.
Let us denote ${\bf Z}_b= (Z_1,\ldots,Z_b)$ for $b=1,2,\ldots,n$, and $T_{lm}= \sum_{j=l+1}^m j Z_j$ for $l=0,1,\ldots,n-1$ and $m=l+1,\ldots,n$.
It follows from the conditioning relation \eqref{CR} that
\begin{equation}
 {\sf P} ({\bf C}^n_b = {\bf a}_b) = {\sf P} ({\bf Z}_b = {\bf a}_b ) \frac{ {\sf P} (T_{bn} = n-a)}{ {\sf P} (T_{0n}=n)}, \label{sta1}
\end{equation}
where $\textbf{a}_b=(a_1,\ldots,a_b)$ and $a= \sum_{j=1}^b j a_j$.

\begin{prop}\label{TDB}
In Case A, if $\theta\geq 1$ for some $n$ and if $\theta^2/n \to 0$, then $ {\sf P} ({\bf C}^n_b = {\bf a}_b) \sim {\sf P} ({\bf Z}_b = {\bf a}_b ) $ for any $\textbf{a}_b$ with any fixed positive integer $b$.
\end{prop}

Proof.
From \eqref{sta1}, in order to prove the desired result, it suffices to show that 
\[\frac{{\sf P}(T_{bn} = n-a)}{{\sf P}(T_{0n}=n)} \to 1.\]
We first calculate $g_{n-a}=\exp(\theta \sum_{j=b+1}^n 1/j) {\sf P}(T_{bn}=n-a)$.
Letting $f(x)=\exp(-\theta \sum_{j=1}^b x^j/j)$, we have
\[ g_{n-a} = [x^{n-a}] (1-x)^{-\theta} f(x) , \]
see equation (5) of \cite{RefABT}.

Let $n$ be a positive integer such that $\theta\geq1$.
It holds that
\begin{eqnarray}
&& [x^{n-a}] (1-x)^{-\theta} f(x)  \nonumber \\ 
& =& [x^{n-a}] (1-x)^{-\theta} \left\{ \sum_{k=0}^{ \lceil \theta \rceil-1} \frac{f^{(k)} (1)}{k!} (-1)^k (1-x)^k + \sum_{k= \lceil \theta \rceil}^\infty \frac{f^{(k)} (1)}{k!} (-1)^k (1-x)^k \right\} \nonumber \\
& =& [x^{n-a}] (1-x)^{-\theta} \left\{ \sum_{k=0}^{\lceil \theta \rceil -1} \frac{f^{(k)} (1)}{k!} (-1)^k (1-x)^k  \right\} + [x^{n-a} ]h(x) ,  \label{DM0}
\end{eqnarray}
where 
\[h(x) = (1-x)^{-\theta} \sum_{k= \lceil \theta \rceil}^\infty \frac{f^{(k)} (1)}{k!} (-1)^k (1-x)^k. \]
Since the right-hand side of \eqref{DM0} is
\begin{equation}
\frac{f(1)(\theta)_{n-a}}{(n-a)!} \left\{ 1 + \sum_{k=1}^{ \lceil \theta \rceil -1} \frac{(-1)^k}{k!} \frac{ f^{(k)}(1) }{f(1)} \frac{(\theta-k)_{n-a}}{(\theta)_{n-a}}   \right\} + [x^{n-a} ]h(x),  \label{DM1}
\end{equation}
the first term and the second term will be evaluated in Lemma \ref{Lem1} and Lemma \ref{Lem2}, respectively.
From Lemma \ref{Lem1}, the elements in the bracket of the first term is $1+O(\theta^2/n)$.
Next we see $[x^{n-a}] h(x)$.
It follows from Lemma \ref{Lem2} that
\[
\left| [x^{n-a}] h(x) \right| 
\leq \frac{1}{r_1^{n}}  \left( \frac{ b e^{(r_2-1)^b} }{r_2} \right)^{\theta} \frac{1}{r_2 -1 -r_1} \{ r_1^a (1+r_1) r_2 \} , 
\]
where $r_1= 1+ c_1 r$ and  $r_2 = 2+ c_2 r$ with constants $c_1,c_2$ such that $1<c_1<c_2$.
By letting $r$ be a positive constant, the right-hand side is $o\left( 1/n^k \right)$ for any positive $k$ since $b$ is fixed and since $\theta^2/n\to0$.
Thence $[x^{n-a}] h(x)=o(1/n)$.

Now we have
\begin{equation}
 g_{n-a} = f(1) \frac{(\theta)_{n-a}}{(n-a)!} \left(1+O \left( \frac{\theta^2}{n}\right) \right) +  o \left( \frac{1}{n} \right)  \label{gna}
\end{equation}
and, as a result,
\begin{eqnarray*}
&& {\sf P} (T_{bn} =n-a) \\
&& = \exp \left(-\theta \sum_{j=1}^n \frac{1}{j} \right)  \frac{ (\theta)_{n-a}}{(n-a)!} \left(  1+ O\left( \frac{ \theta^2}{n} \right) \right)   + \exp \left(-\theta \sum_{j=b+1}^n \frac{1}{j} \right)  o \left( \frac{1}{n} \right) . 
\end{eqnarray*}
On the other hand, 
\begin{equation}
{\sf P} (T_{0n}=n) = \exp\left( -\theta \sum_{j=1}^n \frac{1}{j} \right) [x^n] (1-x)^{-\theta} 
%= \exp \left(-\theta \sum_{j=1}^n \frac{1}{j} \right) \frac{ (\theta)_{n}}{n!} 
= \exp \left(-\theta \sum_{j=1}^n \frac{1}{j} \right) \frac{(\theta)_n}{n!} . \label{t0n}
\end{equation}
If $\theta\to c<\infty$, $(\theta)_n / n!  \sim n^{\theta-1}/\Gamma(\theta)$ and so
\begin{equation}
 \exp\left( \theta \sum_{j=1}^b \frac{1}{j} \right) \frac{n!}{ (\theta)_{n}}
\leq  \exp\left( \theta \log{b} \right) \frac{e^\theta n!}{ (\theta)_{n}} 
\sim  n\Gamma(\theta) \left( \frac{be}{ n} \right)^\theta 
=o(n).  \label{p51t1}
\end{equation}
If $\theta\to\infty$, Lemma \ref{NBL} and the Stirling formula yield that 
\[
\frac{(\theta)_n}{n!}  \sim \frac{n^{\theta-1}}{\Gamma(\theta)} \sim \frac{n^{\theta-1} \theta^{1/2} e^{\theta} }{\sqrt{2\pi} \theta^\theta} , \]
 and hence
\begin{equation}
\exp\left( \theta \sum_{j=1}^b \frac{1}{j} \right) \frac{n!}{ (\theta)_{n}} 
\sim \left( \frac{ \theta e^{\left(\sum_{j=1}^b \frac{1}{j} -1\right)}}{n}  \right)^\theta \frac{ \sqrt{2\pi} n }{  \theta^{1/2} } 
\leq \sqrt{2\pi} \frac{ n  }{  \theta^{1/2} } \left( \frac{ \theta b }{n}  \right)^\theta 
=o(n). \label{p51t2}
\end{equation}
From what has already been proven, we obtain
\begin{eqnarray*}
 \frac{ {\sf P} (T_{bn} =n-a) }{ {\sf P} (T_{0n}=n) } 
&=& \frac{n!}{(\theta)_{n}} \frac{(\theta)_{n-a}}{(n-a)!}\left(  1+ O\left( \frac{\theta^2}{n} \right) \right ) + o(n)   o \left( \frac{1}{n} \right) \\
&\sim& \frac{\theta^\theta}{n^{\theta-1} \theta^{1/2} e^\theta } \frac{(n-a)^{\theta-1} \theta^{1/2} e^\theta }{\theta^\theta} \to 1. 
\end{eqnarray*}
This completes the proof.\\

In the following lemma, we see the first term of \eqref{DM1}.

\begin{lem}\label{Lem1}
Let $f(x)=\exp(-\theta \sum_{j=1}^b x^j/j)$.
For $\theta > 1$, for $k=1,\ldots, \lceil\theta\rceil -1 $ and for any positive integers $a<n$ and $b$, it holds that
\begin{equation}
 \frac{1}{k!} \left| \frac{f^{(k)}(1)}{f(1)} \right| \frac{(\theta - k)_{n-a}}{(\theta)_{n-a}}  \leq \left(\frac{b\theta^2}{n-a}\right)^k. \label{L1G1}
\end{equation}
Moreover, in Case A, if $a=o(n)$, $b=o(n/\theta^2)$, $\theta\geq 1$ for some $n$ and $\theta^2/n\to0$, then
\begin{equation} 
\sum_{k=0}^{\lceil \theta\rceil -1}  \frac{1}{k!} \left| \frac{f^{(k)}(1)}{f(1)} \right| \frac{(\theta - k)_{n-a}}{(\theta)_{n-a}} 
=  1 + \frac{b\theta(\theta-1)}{\theta+n-a} + O\left( \frac{b^2\theta^4}{n^2} \right). \label{L1G2}
\end{equation}
\end{lem}

Proof.
Let $g(x)$ be $-\theta \sum_{j=1}^b x^j/j $.
It holds that
\[ g^{(i)}(x) = -\theta \sum_{j=i}^b \frac{j \cdots (j-i+1)}{j} x^{j-i} = -\theta \sum_{j=1}^{b-i+1} (j)_{i-1} x^{j-1} \]
for $1\leq  i \leq b$.
Thus, for $1\leq i \leq b$, 
\[ 0 \geq g^{(i)}(1) = - \theta  \sum_{j=1}^{b-i+1} (j)_{i-1} \geq -\theta (b-i+1)_{i-1} (b-i+1)  \geq -\theta b^{i} . \]
For  $i > b$, $ g^{(i)}(1) = 0  \geq -\theta b^{i}. $
The Fa\`a di Bruno formula yields that
\[ f^{(k)}(x) = \exp(g(x)) \sum_{j=1}^k B_{k,j}( (g^{(1)} (x), \ldots, g^{(k-j+1)} (x) )) , \]
where $B_{k,j}(\cdot)$ is the partial Bell polynomial, so
\[ \frac{f^{(k)}(1)}{f(1) }  = \sum_{j=1}^k B_{k,j}( (g^{(1)} (1), \ldots, g^{(k-j+1)} (1) )) \]
for any $k=1,2,\ldots$.
By using the triangle inequality,
\begin{eqnarray*}
 | B_{k,j}( (g^{(1)} (1), \ldots, g^{(k-j+1)} (1) )) |
&\leq&  n! \sum_{\{s_\cdot: \sum s_i = j, \ \sum i s_i = k\}} \prod_{i=1}^n \left( \frac{|g^{(i)}(1) |}{i!} \right)^{s_i} \frac{1}{s_i!} \\
&\leq& \theta^j b^{k} \mathcal{S}(k,j),
\end{eqnarray*}
where $\mathcal{S}(k,j)$ is the Stirling number of the second kind.
The above two displays and the triangle inequality imply that
\begin{eqnarray*}
\left| \frac{f^{(k)}(1)}{f(1) } \right| &\leq& \sum_{j=1}^k | B_{k,j}( (g^{(1)} (1), \ldots, g^{(k-j+1)} (1) ))| \\
&\leq& b^{k} \sum_{j=1}^k \theta^j \mathcal{S}(k,j) \\
&\leq& b^{k} \sum_{j=1}^k \theta^j \bar{s}(k,j) = b^{k} (\theta)_k. 
\end{eqnarray*}
For $k \leq \lceil \theta \rceil -1$, the Stirling formula yields that
\[ (\theta)_k \frac{(\theta-k)_{n-a}}{ (\theta)_{n-a}} = \frac{\Gamma(\theta+k)\Gamma(\theta-k +n-a)\Gamma(\theta)}{\Gamma(\theta)\Gamma(\theta-k)\Gamma(\theta +n-a)} \leq \frac{\theta^{2k}}{(n-a)^k}, \]
where we have used $\Gamma(\theta+k)/\Gamma(\theta-k)=(\theta -k)(\theta^2 -(k-1)^2)\cdots (\theta^2 - 1^2) \theta \leq \theta^{2k}$ and $\Gamma(\theta-k+n-a) / \Gamma(\theta+n-a) = 1/ ((\theta-1+n-a)\cdots(\theta-k+n-a)) \leq 1/(n-a)^k $.
We thus have 
\[ \left| \frac{ f^{(k)}(1) }{f(1)} \right| \frac{(\theta-k)_{n-a}}{(\theta)_{n-a}} \leq \left( \frac{b \theta^{2}}{n-a} \right)^k \]
for $k \leq \lceil \theta \rceil -1$, which is \eqref{L1G1}.

Next we prove \eqref{L1G2}.
If $ \theta \leq 1$ for all $n$, the result is obvious because the left-hand side of \eqref{L1G2} is 1.
Otherwise, by letting $n$ be an positive integer such that $\theta > 1$, the desired result follows from
\[ \frac{1}{1!} \frac{f^{(1)}(1)}{f(1)} \frac{(\theta - 1)_{n-a}}{(\theta)_{n-a}} = b\theta \frac{\theta-1}{\theta+n-a}   \]
and from
\begin{equation}
\sum_{k=2}^{\lceil \theta \rceil -1} \frac{1}{k!} \left| \frac{ f^{(k)}(1) }{f(1)} \right| \frac{(\theta-k)_{n-a}}{(\theta)_{n-a}}  
\leq \sum_{k=2}^{\lceil \theta \rceil -1} \left( \frac{b \theta^{2}}{n-a} \right)^k =O\left( \left(  \frac{b \theta^{2} }{n} \right)^2 \right).  \label{T9t}
\end{equation}
This completes the proof.\\

The following lemma is used to evaluate $[x^{n-a}] h(x)$ in \eqref{DM1}.

\begin{lem}\label{Lem2}
Let $f(x)= \exp(-\theta \sum_{j=1}^b x^j/j)$ and let
\[h(x) = (1-x)^{-\theta} \sum_{k= \lceil \theta \rceil}^\infty \frac{f^{(k)} (1)}{k!} (-1)^k (1-x)^k. \]
Then, for any positive integers $a<n$ and $b$, it holds that
\[
\left| [x^{n-a}] h(x) \right| 
\leq \frac{1}{r_1^{n}}  \left( \frac{ b e^{(r_2-1)^b} }{r_2} \right)^{\theta} \frac{1}{r_2 -1 -r_1} \{ r_1^a (1+r_1) r_2 \}
,\]
where $r_1=1+c_1 r$, $r_2=2+c_2 r$, $1< c_1< c_2$, and $r$ is an arbitrary positive constant.
\end{lem}

Proof.
Consider a complex variable ${\bf z} \in \Bbb{C}$.
Since $h({\bf z})$ and $f({\bf z })$ are analytic in $\Bbb{C}$, by using the Cauchy inequality for coefficients, it holds that
\begin{eqnarray*}
 \sup_{|{\bf z}| = r_1} |h({\bf z})| 
&\leq& \sup_{|{\bf z}| = r_1}  \sum_{k=\lceil \theta \rceil}^\infty \left| \frac{f^{(k)}(1) }{k!} \right| |(1-{\bf z})^{k-\theta}|  \nonumber \\
&\leq&  \sum_{k=\lceil \theta \rceil }^\infty \left| \frac{f^{(k)}(1) }{k!} \right|  \sup_{|{\bf z}| = r_1}   |(1-{\bf z})^{k-\theta}|  \nonumber  \\
&\leq&  \sum_{k= \lceil \theta \rceil }^\infty \frac{\sup_{ |1- {\bf z} |  = r_2} |f({\bf z})|}{r_2^k} (1+r_1)^{k-\theta}  \nonumber  \\
&=& \frac{\sup_{ |1- {\bf z} |  = r_2} |f({\bf z})|}{(1+r_1)^{\theta}} \sum_{k= \lceil \theta \rceil }^\infty \left( \frac{ 1+r_1}{r_2}\right)^k  \nonumber  \\
&=&  \frac{ \sup_{ | 1- {\bf z} |  = r_2} |f({\bf z})| }{(1+r_1)^\theta} \frac{ \{ (1+r_1)/r_2 \}^{\lceil \theta \rceil} }{1- (1+r_1)/r_2} .
\end{eqnarray*}
The right-hand side is 
\begin{equation}
\frac{ \sup_{ | 1- {\bf z} |  = r_2} |f({\bf z})| }{ (1+r_1)^{\theta - \lceil \theta \rceil} r_2^{\lceil \theta \rceil} \{ 1 - (1+r_1)/r_2 \} } 
\leq \left( \frac{ b e^{(r_2-1)^b} }{r_2} \right)^{\theta} \frac{(1+r_1) r_2}{r_2 - 1 - r_1} ,  \label{TDa}
\end{equation}
because
\begin{eqnarray*}
\sup_{ |1- {\bf z} |  = r_2} |f({\bf z})| 
&=& \sup_{ |1- {\bf z} |  = r_2} \left|  \exp \left( -\theta \sum_{j=1}^b \frac{{\bf z}^j }{j} \right) \right| \\
&\leq&  \exp \left( \theta \sum_{j=1}^b \frac{( {r_2}-1)^j }{j} \right) \\
&\leq&  \exp \left( \theta \{ \log{b} + (r_2-1)^b \} \right) \\
&=&  (b e^{(r_2-1)^b})^{\theta },
\end{eqnarray*}
where we have used Lemma \ref{SI} for the second inequality.
Hence, it follows from the Cauchy inequality again that
\begin{eqnarray*} 
\left| [x^{n-a}] h(x) \right| 
&\leq& \frac{\sup_{|{\bf z}|=r_1} |h({\bf z})| }{r_1^{n-a}}\nonumber \\
&\leq& \frac{1}{r_1^{n}}  \left( \frac{ b e^{(r_2-1)^b} }{r_2} \right)^{\theta} \frac{1}{r_2 -1 -r_1} \{r_1^a (1+r_1) r_2\} .
\end{eqnarray*}
This completes the proof.\\

Let us provide some remarks on Proposition \ref{TDB}.

\begin{rem}
Proposition \ref{TDB} indicates that when $\theta^2/n \to 0$ the components of $(C_1^n,\ldots,C_b^n)$ are asymptotically independent, and $C_j^n$ asymptotically follows the Poisson distribution with mean $\theta/j$ for $j=1,\ldots,b$.
As a consequence, for any fixed $b$, if $\theta \to c <\infty$ then 
\[
(C_1^n,\ldots,C_b^n) \Rightarrow (Z_1,\ldots,Z_b) ,
\]
and if $\theta\to\infty$ then 
\[
\frac{1}{\sqrt\theta} \left(   C_1^n - \theta , \sqrt{2} \left( C_2^n - \frac{\theta}{2} \right), \ldots, \sqrt{b} \left( C_b^n - \frac{\theta}{b} \right)  \right) \Rightarrow N_b(0,I) ,   
\]
where $N_b(0,I)$ is a $b$-dimensional standard normal variable with independent coordinates.
\end{rem}

\begin{rem}
Proposition \ref{PPAEL} below is stronger than Proposition \ref{TDB}, but the proof is included because some evaluations are different from the proof of Theorem 1 of \cite{RefABT} who used the Darboux lemma (see for instance Theorem of \cite{RefKW}), and because Lemmas \ref{Lem1} and \ref{Lem2} will be used in the proof of Theorem \ref{CPA} below.
\end{rem}

In Proposition \ref{TDB}, $\theta^2/n \to 0$ is assumed.
The following proposition shows that this assumption is necessary for the approximation of $\{ C_j^n \}_{j=1}^b$ by Poisson variables $\{ Z_j \}_{j=1}^b$.

\begin{prop}\label{TDBN}
In Case A, if $\theta \geq 1$ for some $n$, then ${\sf P}({\bf C}^n_b = {\bf a}_b ) \sim {\sf P}({\bf Z}_b = {\bf a}_b )$ for any $\textbf{a}_b$ with any fixed positive integer $b$ only if $\theta^2/n\to 0$ .
\end{prop}

Proof.
To prove the assertion, we see the case that $b=1$.
Let $f(x)=\exp(-\theta x)$, then we have $f^{(k)}(x)=(-\theta)^k f(x)$.
From \eqref{DM1}, $g_{n-a} = [x^{n-a}] (1-x)^{-\theta} f(x)$ equals
\[ \frac{f(1) (\theta)_{n-a}}{(n-a)!} \sum_{k=0}^{\lceil \theta \rceil -1} \frac{\theta^k}{k!} \frac{(\theta-k)_{n-a}}{(\theta)_{n-a}} + [x^{n-a}] h(x). \]
Since
\begin{eqnarray*}
\frac{{\sf P}(T_{1n}=n-a)}{{\sf P}(T_{0n}=n)} 
&=& \frac{ \exp(-\theta \sum_{j=2}^n 1/j) g_{n-a} }{ \exp(-\theta \sum_{j=1}^n 1/j) (\theta)_n/ n!}
= \frac{n! }{(\theta)_n } \frac{ (\theta)_{n-a} }{ (n-a)! } \sum_{k=0}^{\lceil \theta \rceil -1} \frac{\theta^k}{k!} \frac{(\theta-k)_{n-a}}{(\theta)_{n-a}} + o(1) \\
&\sim& \sum_{k=0}^{\lceil \theta \rceil -1} \frac{\theta^k}{k!} \frac{(\theta-k)_{n-a}}{(\theta)_{n-a}}
\end{eqnarray*}
from the proof of Proposition \ref{TDB}, it is enough to show that
\[ \sum_{k=0}^{\lceil \theta \rceil -1} \frac{\theta^k}{k!} \frac{(\theta-k)_{n-a}}{(\theta)_{n-a}} \to 1\]
only if $\theta^2/n \to 0$.

Since $\theta$ is assumed not to decrease as $n$ increases, we study the following three cases: (i) $\theta \geq 2$ for some $n$; (ii) $\theta < 2$ for all $n$ and $\theta>1$ for some $n$; (iii) $\theta \leq 1$ for all $n$.
First, consider (i).
Let $n$ be a positive integer such that $\theta \geq 2$.
Then, it holds that
\begin{eqnarray*}
\sum_{k=0}^{\lceil \theta \rceil -1} \frac{\theta^k}{k!} \frac{(\theta-k)_{n-a}}{(\theta)_{n-a}} 
&\geq& \sum_{k=0}^{\lceil \theta \rceil -2} \frac{\theta^k}{k!} \frac{(\theta-k)_{n-a}}{(\theta)_{n-a}}  \\
&\geq& \sum_{k=0}^{\lceil \theta \rceil -2} \frac{\theta^k}{k!} \frac{(\lceil \theta \rceil -1 -k)_{n-a}}{(\lceil \theta \rceil -1)_{n-a}} \\
&=&  \sum_{k=0}^{\lceil \theta \rceil -2} \frac{\theta^k}{k!} \frac{(\lceil \theta \rceil -2)!}{ (\lceil \theta \rceil  -k -2)! }  \frac{(\lceil \theta \rceil -2 -k + n -a)!}{(\lceil \theta \rceil -2 +n-a)!} \\
&=&  \sum_{k=0}^{\lceil \theta \rceil -2} \binom{\lceil \theta \rceil -2}{k} \theta^k \frac{ 1}{ (\lceil \theta \rceil - 2 -k +1  + n-a ) \cdots  (\lceil \theta \rceil - 2   + n-a ) } \\
&=&  \sum_{k=0}^{\lceil \theta \rceil -2} \binom{\lceil \theta \rceil -2}{k} \left( \frac{\theta}{ n } \right)^k \frac{ 1}{ \left(1 + \frac{\lceil \theta \rceil - 2 -k +1  -a}{n} \right) \cdots  \left(1 + \frac{\lceil \theta \rceil - 2   -a}{n} \right) } \\
&\geq&  \sum_{k=0}^{\lceil \theta \rceil -2} \binom{\lceil \theta \rceil -2}{k} \left\{ \frac{\theta}{ n \left(1 + \frac{\lceil \theta \rceil - 2   -a}{n} \right) }  \right\} ^k,
\end{eqnarray*}
where we have used Lemma \ref{SIL} for the second inequality.
From the binomial theorem, the right-hand side is equal to
\[
\left\{  1 + \frac{\theta}{ n \left( 1 + \frac{\lceil \theta \rceil - 2   -a}{n} \right) } \right\}^{\lceil \theta \rceil - 2}
=\left[ \left\{  1 + \frac{\theta}{ n \left( 1 + \frac{\lceil \theta \rceil - 2   -a}{n} \right) } \right\}^{n/\theta} \right]^{ \theta (\lceil \theta \rceil - 2)/n}.
\]
The above display is not less than 1 and converges to 1 only if $\theta^2/n \to 0$.
Second, consider (ii).
Let $n$ be a positive integer such that $\theta > 1$.
Then, it holds that
\[ \sum_{k=0}^{\lceil \theta \rceil -1} \frac{\theta^k}{k!} \frac{(\theta-k)_{n-a}}{(\theta)_{n-a}}  = 1 + \frac{\theta(\theta -1)}{ \theta + n -a -1 }, \]
which converges to 1 only if $\theta^2/n \to 0$.
Finally, consider (iii).
Let $n$ be a positive integer such that $\theta = 1$.
Then, it holds that
\[ \sum_{k=0}^{\lceil \theta \rceil -1} \frac{\theta^k}{k!} \frac{(\theta-k)_{n-a}}{(\theta)_{n-a}}  = 1 . \] 
This completes the proof.\\

Thence, we have the following corollary to Propositions \ref{TDB} and \ref{TDBN}.

\begin{cor}
In Case A, if $\theta \geq 1$ for some $n$, then ${\sf P}({\bf C}^n_b = {\bf a}_b) \sim {\sf P}({\bf Z}_b = {\bf a}_b )$ for any $\textbf{a}_b$ with any fixed positive integer $b$ if and only if $\theta^2/n\to 0$.
\end{cor}

Subsequently, let us derive the result corresponding to \eqref{AST} following a similar programme to \cite{RefAST}.
It follows from \eqref{CR} that
\[
d_b(n) = \sum_{a =0}^\infty {\sf P}(T_{0b}=a ) \left( 1- \frac{ {\sf P} (T_{bn}=n-a)}{ {\sf P} (T_{0n}=n)}  \right)^+  , 
\]
see (50) of \cite{RefAST}.
Firstly, via the large deviation inequality, we see that $d_b(n)$ is approximated by
\[  
 \sum_{a=0}^{ \lfloor J_n \rfloor}  {\sf P}(T_{0b}=a ) \left( 1- \frac{ {\sf P} (T_{bn}=n-a)}{ {\sf P} (T_{0n}=n)}  \right)^+ 
\]
with 
\[J_n = \min(b \theta \log{n}, b^{2/3} (\theta n )^{1/3}). \]
From the definition, if $1\leq b \leq n \theta^{-2}(\log{n})^{-3}$ then $J_n=b \theta \log{n}$ and otherwise $J_n=b^{2/3} (\theta n )^{1/3} = b (\theta n/b)^{1/3}$.
In contrast to \cite{RefAST}, $J_n$ includes $\theta$ since we consider $\theta \to \infty$, but a similar treatment perform well.

\begin{lem}\label{ASTL9}
In Case A, with $b=o(n/\theta^2)$, it holds that
\[
\sum_{a > J_n} {\sf P}(T_{0b}=a ) \left( 1- \frac{ {\sf P} (T_{bn}=n-a)}{ {\sf P} (T_{0n}=n)}  \right)^+ \leq  {\sf P}(T_{0b}> J_n) = o\left( \left( \frac{b}{n}  \right)^k \right)
\]
for any positive $k$.
\end{lem}

Proof.
The first inequality is obvious, so we see the latter one.
From Lemma 8 of \cite{RefAST}, for any $b\geq1,w>0$, it holds that
\begin{equation}
\log {\sf P}(T_{0b} \geq bw) \leq \log (\theta e/w)^w. \label{LDp1}
\end{equation}
If $1\leq b \leq n \theta^{-2}(\log{n})^{-3}$ then, by putting $w=\theta\log{n}$, the right-hand side of \eqref{LDp1} is
\[ (\theta \log{n}) (1-\log\log{n}) \sim - \theta (\log\log{n}) \log{n}  \]
which tends to minus infinity faster than $-k \log{n}$ for any positive $k$.
If $b \geq n \theta^{-2}(\log{n})^{-3}$ then, by putting $w=(\theta n/b)^{1/3}$, the right-hand side of \eqref{LDp1} is
\[ \left( \frac{\theta n }{b} \right)^{1/3} \left( 1-\frac{1}{3} \log\left(\frac{n}{b} \right) \right) \sim - \frac{1}{3} \left(\frac{\theta n }{b}\right)^{1/3} \log\left(\frac{n}{b} \right)   \]
which tends to minus infinity faster than $-k \log{(n/b)}$ for any positive $k$.
This completes the proof.\\

The next lemma shows that $ (|1-\theta|/n) {\sf E}\left[  (T_{0b} - \theta b )^+ 1\{ T_{0b} \leq J_n  \}  \right]$ is approximately $(|1-\theta|/n) {\sf E}\left[  ( T_{0b} - \theta b )^+ \right] $.

\begin{lem}\label{ASTL10}
In Case A, if $\theta^2/n \to 0$ then it holds that
\[ \frac{|1-\theta|}{n} {\sf E}\left[  |T_{0b} - \theta b | 1\{ T_{0b} > J_n  \}  \right]  = o\left( \left( \frac{b}{n}  \right)^k \right) \]
for $b=o(n/\theta^2)$ and for any positive $k$.
\end{lem}

Proof.
From the Schwartz inequality, it follows that
\begin{eqnarray*}
&& \frac{|1-\theta|}{n} {\sf E}\left[  |T_{0b} - b \theta | 1\{ T_{0b} > J_n  \}  \right]  \\
&\leq& \left( \frac{|1-\theta|^2}{n^2} {\sf E}[(T_{0b} - \theta b)^2] {\sf P}(T_{0b}>J_n) \right)^{1/2} \\
&\leq& \left( \frac{|1-\theta|^2}{n^2} \theta b^2 {\sf P}(T_{0b}>J_n) \right)^{1/2} \\
&\leq& \left[ \left\{  \frac{(1+\theta)^2 b  }{n} \right\}^2 {\sf P}(T_{0b}>J_n) \right]^{1/2} 
\end{eqnarray*}
where we have used ${\sf E}[(T_{0b} - {\sf E}[T_{0b}])^2] ={\rm var}(T_{0b}) = \sum_{j=1}^b j^2 (\theta/j) = \theta \sum_{j=1}^b j \leq \theta b^2$ for the second inequality.
Lemma \ref{ASTL9} yields that ${\sf P}(T_{0b}>J_n) = o((b/n)^{2k})$ for any positive $k$.
This completes the proof.\\

The following result is an extension of \eqref{AST} to large $\theta$ setup. 

\begin{thm}\label{CPA}
In Case A, if $\theta \geq 1$ for some $n$ and $\theta^2/n \to 0$, then 
\begin{equation}
 d_b(n) = \frac{\theta-1}{2n} {\sf E}[| T_{0b} -\theta b  |]  + o\left( \frac{b\theta^2}{n} \right) \label{TAR}
\end{equation}
for 
\begin{equation}
b = o \left( \frac{n}{\theta^2 \log{n} } \right). \label{bass}
\end{equation}
In addition, when $\theta\to\infty$, it holds that $d_b(n)=o(b\theta^2/n)$.
\end{thm}

Proof.
Let $n$ be a positive integer such that $\theta \geq 1$.
Since it follows from Lemma \ref{ASTL9} that
\[ d_b(n)=  \sum_{a=0}^{ \lfloor J_n \rfloor}  {\sf P}(T_{0b}=a ) \left( 1- \frac{ {\sf P} (T_{bn}=n-a)}{ {\sf P} (T_{0n}=n)}  \right)^+     + o \left( \left( \frac{b}{n} \right)^k \right) \]
for any positive $k$, we see the first term.

Let $g_{n-a} = \exp(\theta \sum_{j=b+1}^n 1/j) {\sf P}(T_{bn} =n-a)$.
For $a \leq \lfloor J_n \rfloor$, we have
\begin{eqnarray*}
g_{n-a} 
&=& \frac{f(1)(\theta)_{n-a}}{(n-a)!} \left\{ 1 + \frac{b\theta (\theta-1)}{\theta+n-a} + O\left( \frac{b^2 \theta^4}{n^2} \right)  \right\}  +   [x^{n-a}]h(x) \\
&=&  \frac{f(1)(\theta)_{n-a}}{(n-a)!} \left\{ 1 + \frac{b\theta (\theta-1)}{n}\left( 1 +  O \left( \frac{a+\theta}{n}\right)  \right) + O\left( \frac{b^2 \theta^4}{n^2} \right)  \right\}  +   [x^{n-a}]h(x)
\end{eqnarray*}
and the last term should be evaluated for $a$ and $b$ growing with $n$.

If $b$ does not diverge, as it is seen in the proof of Proposition \ref{TDB}, $ [x^{n-a}]h(x) = o(1/n^2)$ since $a/n \leq J_n/n \to 0$.
Thence, we consider the case that $b\to\infty$.
Using Lemma \ref{Lem2} with $r = 1/b$, we have
\[
\left| [x^{n-a}] h(x) \right| 
\leq \frac{b}{(1+c_1/ b)^{n-a}}  \left( \frac{ b e^{(1+c_2/b)^b} }{2+c_2/ b} \right)^{\theta} \frac{ (2+c_1/ b)(2+c_2/ b) }{(c_2-c_1)}
. \]
Since $(1+c_1/ b)^{n-a} \sim e^{(n-a)c_1/ b}$ and since $(1+c_2/b)^b \sim e^{c_2}$, the right hand side is asymptotically equal to
\begin{eqnarray*}
 \frac{b^{\theta+1} A_1^\theta }{ \exp({(n-a)c_1/ b}) } A_2 
 &=& A_2 \exp \left( (\theta+1)\log{b} + \theta \log{A_1} -  \frac{(n-a)c_1}{b} \right) \nonumber \\ 
 &=& A_2 \exp \left( \frac{n}{b} \left\{ \frac{b(\theta+1)\log{n}}{n} \frac{\log{b}}{\log{n}} + \frac{b\theta}{n} \log{A_1} -  \left(  1- \frac{a}{n} \right)c_1 \right\} \right),  
 \end{eqnarray*}
where $A_1= (\exp({e^{c_2}}))/2$ and $A_2 = 4/(c_2-c_1)$.
From \eqref{bass}, the right-hand side is
\begin{equation}
A_2 \exp \left( - \frac{n}{b} \left( c_1 + o(1) \right) \right) = A_2 \exp \left( - \frac{n\left( c_1 + o(1) \right)}{b\log{n}} \log{n}  \right) = A_2 n^{-\frac{n\left( c_1 + o(1) \right)}{b\log{n}}} ,  \label{rel}
\end{equation}
where we have used $a/n \leq J_n/n \to 0$.
The right-hand side is $ o \left( 1/n^k \right)$ for any positive constant $k$.
After all, we have $\left| [x^{n-a}] h(x) \right|  = o \left( b/n^2 \right)$ even when $b\to\infty$.

Now we have
\[g_{n-a} =  \frac{f(1)(\theta)_{n-a}}{(n-a)!} \left\{ 1 + \frac{b\theta (\theta-1)}{n} \left( 1 +  O \left( \frac{a+\theta}{n}\right)  \right) + O\left( \frac{b^2 \theta^4}{n^2} \right)  \right\}  +  o\left( \frac{b}{n^2} \right) . \]
This expansion, $f(1) = \exp(-\sum_{j=1}^b 1/j)$, ${\sf P}(T_{bn} =n-a) = \exp \left(-\theta \sum_{j=b+1}^n 1/j \right)g_{n-a} $ and \eqref{t0n}--\eqref{p51t2} yield that
\[ \frac{ {\sf P} (T_{bn} =n-a) }{ {\sf P} (T_{0n}=n) } 
= \frac{n!(\theta)_{n-a}}{(n-a)!(\theta)_{n}}  \left\{  1+ \frac{b\theta (\theta-1)}{n}  + o\left( \frac{b \theta^2}{n} \right) \right\} + o(n)   o \left( \frac{b}{n^2} \right) . \]
Since $a\theta/n \leq  aJ_n/n \to 0$ and $a^2/(nb) \leq J_n^2/(nb) \to 0$ which follow from $\theta/J_n \to 0$ and \eqref{bass}, the binomial expansion and Lemma \ref{NBL} yield that
\begin{eqnarray*}
\frac{ (\theta)_{n-a}/(n-a)!  }{(\theta)_n/n!} 
&=& \frac{(n-a)^{\theta-1} \left\{ 1+ \frac{\theta(\theta - 1)}{2(n-a)} + O \left( \frac{\theta^4}{n^2} \right) \right\} /\Gamma(\theta) }{ n^{\theta-1} \left\{ 1 + \frac{\theta(\theta - 1)}{2n} +O \left( \frac{\theta^4}{n^2} \right)  \right\} /\Gamma(\theta) } \\
&=& \left( 1- \frac{a}{n} \right)^{\theta-1} \frac{ \left\{ 1+ \frac{\theta(\theta - 1)}{2n} + O \left( \frac{\theta^2(a+\theta^2)}{n^2} \right) \right\} }{ \left\{ 1 + \frac{\theta(\theta - 1)}{2n} +O \left( \frac{\theta^4}{n^2} \right)  \right\}  } \\
&=& \left\{ 1- \frac{a(\theta-1)}{n} + O\left( \frac{a^2\theta^2}{n^2} \right) \right\} \left( 1+O \left( \frac{\theta^2(a+\theta^2))}{n^2} \right)  \right),
\end{eqnarray*}
and hence
\[ \frac{n!(\theta)_{n-a}}{(n-a)!(\theta)_{n}} = 1 - \frac{(\theta-1)a}{n} + O\left( \frac{\theta^2 (a^2+\theta^2)}{n^2} \right)  =1 - \frac{(\theta-1)a}{n} + o \left( \frac{b \theta^2}{n} \right). \]
Therefore, it holds that
\begin{eqnarray*}
\frac{ {\sf P} (T_{bn} =n-a) }{ {\sf P} (T_{0n}=n) } &=&  
\left\{  1 - \frac{(\theta-1)a}{n} + o \left( \frac{b \theta^2}{n} \right) \right\}  \left\{  1+ \frac{b\theta (\theta-1)}{n}  + o\left( \frac{b \theta^2}{n} \right) \right\} + o \left( \frac{b}{n} \right) \\
&=& 1-  \left\{  \frac{(\theta-1)a}{n}   - \frac{b\theta (\theta-1)}{n} \right\} + o\left( \frac{b \theta^2}{n} \right).
\end{eqnarray*}
From what has already been proven, it holds that
\begin{eqnarray*}
d_b(n) &=& \sum_{a=0}^{ \lfloor J_n \rfloor}  {\sf P}(T_{0b}=a )\left( \frac{(\theta-1)a}{n}   - \frac{b\theta (\theta-1)}{n}  \right)^+    + o\left( \frac{b \theta^2}{n} \right) \\
&=& \frac{1}{n} \sum_{a=0}^{ \infty }  {\sf P}(T_{0b}=a )\left((\theta-1)( a-b\theta ) \right)^+ 1\{ a \leq J_n \}    + o\left( \frac{b \theta^2}{n} \right) \\
&=&  \frac{1}{n}  {\sf E} \left[\left((\theta-1) (T_{0b} -b\theta ) \right)^+ 1\{ T_{0b} \leq J_n \}  \right]  + o\left( \frac{b \theta^2}{n} \right) \\
&=&  \frac{1}{n}  {\sf E} \left[\left((\theta-1) (T_{0b} -b\theta ) \right)^+ \right]  + o\left( \frac{b \theta^2}{n} \right) \\
&=&  \frac{\theta-1}{2n}  {\sf E} \left[\left| T_{0b} -b\theta \right| \right]  + o\left( \frac{b \theta^2}{n} \right) 
\end{eqnarray*}
where we have used Lemma \ref{ASTL10} in the fourth equality and the relation ${\sf E}[(T_{0b} - b\theta)^+] = {\sf E} [|T_{0b} - b\theta|]/2$, which follows from $ {\sf E} [T_{0b} - b\theta]=0$, in the fifth equality.

Finally, consider the case that $\theta\to\infty$.
It follows from the Jensen inequality that
\[ {\sf E} \left[\left| T_{0b} -b\theta  \right| \right] \leq \sqrt{ {\sf E} \left[\left| T_{0b} -b\theta  \right|^2 \right]} = O(\theta^{1/2} b), \]
which implies $d_b(n)=o(b\theta^2/n)$.
This completes the proof.\\

Next we discuss the Poisson process approximation via the Feller coupling (see Subsection \ref{sec22}).
The following result follows directly from \eqref{wb1}.

\begin{prop}\label{PPAEL}
Suppose that $\theta > 1$ for all $n$ and $\theta^2/n \to 0$.
In Case A, $d_b^W(n) \to 0$ if, and only if, $b=o(n/\theta^2)$.
In addition, when $\theta\to\infty$, it holds that 
\[ d_b^W(n) \sim \frac{\theta^2}{n} b. \]
\end{prop}

\begin{rem}
When $\theta\to\infty$, Theorem \ref{CPA} and Proposition \ref{PPAEL} yield that
\[ d_b(n)=o \left( \frac{b\theta^2}{n} \right), \quad d_b^W(n)=\Theta \left( \frac{b\theta^2}{n} \right) \]
with $b=o(n/(\theta^2\log{n}))$, which shows that the asymptotic decay rates of $d_b(n)$ and $d_b^W(n)$ are different.
\end{rem}

Let $S^k_n$ be the $k$-th shortest cycle length in a Ewens partition, that is,
\[
S^k_n = \inf(j: C_1^n + \ldots+C_j^n \geq k) 
\]
for $k=1,2,\ldots$ and $S^k_n = \infty$ when there is no such $j$.
See Section 2E of \cite{RefAT}.
The following statement is a direct corollary to Theorem \ref{CPA}.

\begin{cor}\label{SCA}
Let $r$ be a positive integer such that $r = o(n/\theta^2)$ and let $\delta_r = \sum_{j=1}^r \theta/j$.
Under the assumption of Proposition \ref{PPAEL}, 
\[  {\sf P} (S^k_n \leq r )  \sim \sum_{x=0}^{k-1} e^{-\delta_r} \frac{\delta_r^x}{x!}. \]
\end{cor}

Proof.
Proposition \ref{PPAEL} yields that
\[  {\sf P} (S^k_n \leq r ) = {\sf P} \left( \sum_{j=1}^r C_j^n < k \right)  \sim {\sf P} \left(  \sum_{j=1}^r Z_j < k \right)  = \sum_{x=0}^{k-1} e^{-\delta_r} \frac{\delta_r^x}{x!} . \]
This completes the proof.\\

\begin{rem}
Corollary \ref{SCA} yields that, under the assumption of Proposition \ref{PPAEL}, ${\sf P}(S_n^1 =1) \sim e^{-\theta}$, so if $\theta\to \infty$ then ${\sf P}(S_n^1 =1) \to 1$ and  if $\theta \to c <\infty$ then ${\sf P}(S_n^1 =1) \to e^{-c}<1$.
Note that when the Pitman sampling formula is considered, the shortest cycle length converges to 1 in probability except the Ewens sampling formula (see \cite{RefM}).
\end{rem}

The uniform bound with respect to $b$, which gives an extension of \eqref{ub}, is given in the following proposition.
Its applications to functional central limit theorems will be presented in the next section.

\begin{prop}\label{PPAU}
In Case A,
\[ d_n^W(n)  = O(\theta \log{(1+\theta)}). \]
\end{prop}

Proof.
By using the triangle inequality and \eqref{wbA}, it holds that
\begin{eqnarray}
d_n^W(n)  &\leq &  \sum_{j=1}^n {\sf E} \left[ |C_j^n - C_j^\infty | \right]  \\
&\leq& \sum_{j=1}^b {\sf E} \left[ |C_j^n  - C_j^\infty| \right] + \sum_{j=b+1}^n {\sf E} \left[ C_j^n \right] +  \sum_{j=b+1}^n {\sf E} \left[ C^\infty_j \right] \nonumber \\
&\leq& \frac{b\theta(\theta+1)}{\theta +n - b} + 1 +2\theta\log{\left( \frac{n}{b}\right)}, \label{ubp1}
\end{eqnarray}
for any $b=1,2,\ldots,n$, see the proof of Theorem 2 of \cite{RefABT}.
When $\theta \to \infty$, by setting $b = \lfloor n/\theta \rfloor$, the first and third terms in \eqref{ubp1} are $O(\theta)$ and $O(\theta \log{\theta})$, respectively.
Otherwise, by setting $b=\lfloor n/2 \rfloor$ the result holds with the bound $d_n^W(n)=O(1)$.
Hence, $d_n^W(n) =  O((\theta \log{\theta}) \vee 1) = O(\theta \log{(1+\theta)}).$
This completes the proof.\\

\subsection{Case C}
The probability mass function \eqref{ESF} is obtained from the conditioning relation \eqref{CR} with a sequence of independent Poisson variables with respective means $\theta/j$.
We also get \eqref{ESF} from \eqref{CR} with Poisson varivables with respective means $(\theta/j) (n/\theta)^j$ (see for instance \cite{RefW1}).
In Case C, the following lemma shows that ${\sf E}[Z_j]=(\theta/j) (n/\theta)^j$ rather fit.

\begin{lem} \label{CM}
In Case C,
\[ {\sf E} \left[ C_j^n \right] \sim \frac{\theta}{j} \left( \frac{n}{\theta} \right)^j. \]
Therefore,  for $j=2,3,\ldots$, if $\theta(n/ \theta)^j \to0$, then
\begin{equation}
C_j^n \to^p 0.  \label{CM2} 
\end{equation}
\end{lem}

Proof.
It holds that
\[ {\sf E}[C_j^n] = \frac{\theta}{j} \frac{n! }{(n-j)! } \frac{\Gamma(n+\theta-j)}{\Gamma(n+\theta)} \]
which is (2.18) of \cite{RefW1}.
Since the Stirling formula $\Gamma(x) = \sqrt{2\pi} x^{x-1/2} e^{-x} +O(x^{x-3/2} /e^{x})$ as $x\to\infty$ yields that
$ \Gamma(x-c)/\Gamma(x) \sim x^{-c} $
as $x\to\infty$ for any $c<x$, it holds that
\[
{\sf E}[C_j^n] \sim \frac{\theta}{j} n^j \frac{1}{(n+\theta)^j} \sim \frac{\theta}{j} \left( \frac{n}{\theta} \right)^j. \]
Hence, the result \eqref{CM2} follows from $C_j^n \geq 0$.
This completes the proof.

\begin{rem}
Since
\[ \frac{n! }{(n-j)! n^j} \leq \frac{n! }{(n-j)! (n-j+1)_{j}} = 1 \]
and
\[ \frac{\Gamma(n-j+\theta) \theta^j}{\Gamma(n+\theta)}  = \frac{\theta^j}{(n-j + \theta)_j} \leq 1, \]
it holds that 
\begin{equation}
 \frac{{\sf E}[j C_j^n ]}{\theta \left( n/\theta \right)^j } = \frac{n! }{(n-j)! n^j} \frac{\Gamma(n-j+\theta) \theta^j}{\Gamma(n+\theta)} \leq 1  \label{mineq}
\end{equation}
for $j=1,2,\ldots,n$.
\end{rem}

According to Lemma \ref{CM}, it may be natural to consider that $C_j^n$ and Poisson variable with mean $(\theta/j) (n/\theta)^j$ are asymptotically similar, but Proposition \ref{P3B} indicates that, except Case C3, an independent process approximation by Poisson variables with means $(\theta/j) (n/\theta)^j$ seems difficult in the sense of the joint distribution.
Actually, the following theorem shows that in Case C2 the linear relation $n-(C_1^n +2C_2^n)\Rightarrow 0$ between $C_1^n$ and $C_2^n$ asymptotically remains.

\begin{thm}\label{PPAC}
(i) In Case C2, 
\[ (  C_1^n - n , C_2^n) \Rightarrow (-2 P_{c/2}, P_{c/2}), \]
where $P_{c/2}$ is a Poisson variable with ${\sf E}[P_{c/2}] = \lim_{n,\theta} \{ n^2/(2\theta) \} = c/2$, and $\sum_{j=3}^n | C_j^n | \Rightarrow 0$.\\
(ii) In Case C3, 
\[\sum_{j=1}^n |C_j^n - n 1\{j=1\}| \Rightarrow 0.\]
\end{thm}

Proof.
(i) It follows from \eqref{mineq} that
\[ {\sf E} \left[ \sum_{j=3}^n jC_j^n \right] \leq \sum_{j=3}^n \theta \left( \frac{n}{\theta} \right)^j =  \frac{n^3}{\theta^2} \frac{ 1- (n/\theta)^{n-2} }{1-n/\theta}, \]
which implies that
\[ n - (C_1^n +2C_2^n) = \sum_{j=3}^n j C_j^n  \to^p 0. \]
It yields that $\sum_{j=3}^n |C_j^n| \to^p 0$ and so 
\[K_n - (C_1^n +C_2^n ) \to^p 0 . \]
These displays yield that 
\[ K_n -n + C_2^n \to^p 0. \]
From this, we obtain
\begin{eqnarray*}
&& K_n -n - \left( \theta\log\left(1+\frac{n}{\theta}\right) -n \right) +  C_2^n +  \left( \theta\log\left(1+\frac{n}{\theta}\right) -n \right) \\
&=& K_n  -\theta\log\left(1+\frac{n}{\theta}\right)  +  C_2^n - \left(n- \theta\log\left(1+\frac{n}{\theta}\right) \right) \\
&\to^p& 0.
\end{eqnarray*}
We conclude from \eqref{TC}, which means $K_n  -\theta\log(1+n/\theta)  \Rightarrow c/2 - P_{c/2}$, that
\[  C_2^n - \left(n- \theta\log\left(1+\frac{n}{\theta}\right) \right) = C_2^n - \frac{n^2}{2\theta} + o(1)  \Rightarrow P_{c/2} - \frac{c}{2},  \]
hence that $C_2^n \Rightarrow P_{c/2}$.
Moreover, from what has already been proven, we obtain $n-C_1^n \Rightarrow 2P_{c/2}$.

(ii) Since
\[ {\sf E}[|n-C_1^n |] = {\sf E}[n-C_1^n] = n- \frac{n }{1+ (n-1)/\theta}  =  \frac{n^2}{\theta} + O\left( \frac{n^3}{\theta^2}  \right), \]
Lemma \ref{CM} yields the result.
This completes the proof.\\

Theorem \ref{PPAC} directly implies the following corollaries which represent properties of the shortest cycle length $S_n=S_n^1$ and the longest cycle length $L_n$ in a Ewens partition. 
These extreme sizes are of interest in the combinatorial context, see for instance \cite{RefM}.

\begin{cor}
In Case C2 or C3,
\[ {\sf P} (S_n =1) \to 1. \]
\end{cor}

Proof.
In Case C3, the conclusion is obvious, so we see Case C2.
It holds that
\begin{eqnarray*}
{\sf P}(S_n\geq 2) &=& {\sf P} (C_1^n =0) ={\sf P} (n-C_1^n \geq n) \\
&\sim& {\sf P} (2P_{c/2} \geq n) = 1 - \sum_{k < n/2} e^{- \frac{c}{2}} \left(\frac{c}{2} \right) ^k \frac{1}{k!} \\
&\to& 0. 
\end{eqnarray*}
This completes the proof.

\begin{cor}\label{LAR}
In Case C2, for a positive integer $r$
\[ {\sf P}(L_n \leq r) \sim 
\begin{cases} 
e^{-c/2} & (r=1) \\
1. & (r\geq 2)
\end{cases}
\]
\end{cor}

Proof.
For $r=1$, it holds that 
\[  {\sf P}(L_n \leq 1)  = {\sf P}(C_1^n = n) \sim {\sf P}(W=0) =  e^{-c/2}.  \]
For $r\geq 2$, it holds that $ {\sf P}(L_n \leq r)  \sim 1. $
This completes the proof.

\begin{rem}
Since the law of the singleton $C_1^n$ in a Ewens partition is given by
\[ {\sf P}(C_1^n=k) = \frac{\theta^k}{k!} \left\{ \sum_{j=0}^{n-k} (-1)^j \frac{\theta^j}{j!} \frac{(n+1-k-j)_{k+j}}{(n+\theta-k-j)_{k+j}} \right\} \]
for $k=0,1,\ldots,n$, Corollary \ref{LAR} directly follows from
\begin{eqnarray*}
{\sf P}(C_1^n = n) &=& \frac{\theta^n}{(\theta)_n} = 1 \frac{1}{1+1/\theta} \cdots \frac{1}{1+(n-1)/\theta} = 1 -\frac{1}{\theta} \frac{n(n-1)}{2} +o(1) \\
&=& e^{-n^2/(2\theta)} + o(1) . 
\end{eqnarray*}
\end{rem}

%%%%%%%%%%%%%%%%%%%%%%%%%%%%%%%%%%%%
\section{Functional central limit theorems}\label{sec6}
As a corollary to Proposition \ref{PPAU}, the functional central limit theorems which extend the results of \cite{RefH} and \cite{RefT1} follow.
Before the results, as a corollary to Proposition \ref{PPAU}, let us states the error bounds of Poisson process approximations in the sense of the expectation of the error in the supremum norm and in the $L^2$ norm.

\begin{cor}\label{EBA}
In Case A, 
\begin{equation}
 {\sf E}\left[ \sup_{u \in [0,1]} \left| \sum_{j=1}^{\lfloor n^u \rfloor} \frac{(C_j^n - C_j^\infty)}{\sqrt{\theta\log{n}}} \right| \right] =O \left( \sqrt{ \frac{\theta}{\log{n}}}  \log{(1+\theta)}  \right) \label{EB1}
\end{equation}
and
\begin{equation}
{\sf E}\left[ \left[ { \int_0^1 \left\{ \frac{\sum_{j=1}^{\lfloor n^u \rfloor} (C_j^n - C_j^\infty)}{\sqrt{\sum_{j=1}^{\lfloor n^u \rfloor} \theta/j } } \right\}^2 du } \right]^{1/2} \right] = O \left(  \sqrt{\frac{\theta \log\log{n} }{\log{n}}}  \log{(1+\theta)} \right) .  \label{EB2}
\end{equation}
\end{cor}

Proof.
The result \eqref{EB1} follows from 
\[ {\sf E}\left[ \sup_{u \in [0,1]} \left| \sum_{j=1}^{\lfloor n^u \rfloor} \frac{(C_j^n - C_j^\infty)}{\sqrt{\theta\log{n}}} \right| \right] 
\leq {\sf E}\left[  \sup_{u \in [0,1]} \sum_{j=1}^{\lfloor n^u \rfloor} \frac{ \left| C_j^n - C_j^\infty \right|  }{\sqrt{\theta\log{n}}} \right]
=  \frac{ \sum_{j=1}^{ n } {\sf E} \left[ \left| C_j^n - C_j^\infty \right| \right] }{\sqrt{\theta\log{n}}} ,  \]
and Proposition \ref{PPAU}.
The result \eqref{EB2} follows from
\[ \int_0^1 \left|  \frac{\sum_{j=1}^{\lfloor n^u \rfloor}(C_j^n - C_j^\infty)}{\sqrt{\sum_{j=1}^{\lfloor n^u \rfloor} \theta/j }} \right|^2 du 
\leq \frac{2\left(\frac{1}{\log2} - \log\log{2} +\log\log{n} \right) }{\theta\log{n}}  \left( \sum_{j=1}^n |C_j^n - C_j^\infty| \right)^2 \]
(for this evaluation see the proof of Lemma 3.1 of \cite{RefT1}) and from Proposition \ref{PPAU}.
This completes the proof.\\

By using Corollary \ref{EBA}, we provide functional central limit theorems which slightly extend the preceding result in which $\theta$ is assumed to be fixed.

\begin{prop} \label{FCLT1}
(i) In Case A, if 
\begin{equation}
\frac{\theta}{\log{n}} \left(\log{\theta}  \right)^2 \to 0, \label{FASS} 
\end{equation}
 then the random process $X^1_n(\cdot)$ defined in \eqref{RF0} converges weakly to a standard Brownian motion $(B(u))_{0\leq u \leq 1}$ in $D[0,1]$.\\
(ii) In Case A, if 
\begin{equation}
\frac{\theta \log\log{n}}{\log{n}} \left(\log{\theta}  \right)^2 \to 0, \label{FASS1} 
\end{equation}
then both of the random processes $X^2_n(\cdot)$ and $X^3_n(\cdot)$, which are defined in \eqref{RF1} and \eqref{RF2} respectively, converge weakly to $(B(u)/\sqrt{u})_{0<u<1}$ in $L^2(0,1)$.
\end{prop}

Proof.
(i)  From \eqref{EB1} and the assumption \eqref{FASS}, it follows that
\begin{equation}
\sup_{u \in [0,1]} \left| \sum_{j=1}^{\lfloor n^u \rfloor} \frac{(C_j^n - C_j^\infty)}{\sqrt{\theta\log{n}}} \right| 
\leq \sup_{u \in [0,1]} \sum_{j=1}^{\lfloor n^u \rfloor} \frac{ | C_j^n - C_j^\infty | }{\sqrt{\theta\log{n}}}
= \sum_{j=1}^{n} \frac{ | C_j^n - C_j^\infty | }{\sqrt{\theta\log{n}}}
\to^p 0. \label{FC1}
\end{equation}
By using the functional central limit theorem for Poisson processes in $D[0,1]$, the random process
\[ \left( \frac{ \sum_{j=1}^{\lfloor n^u \rfloor} C_j^\infty - \sum_{j=1}^{\lfloor n^u \rfloor} \theta/j  }{ \sqrt{ \sum_{j =1}^{n} \theta/j  } } \right)_{0\leq u\leq1} \]
converges weakly to a standard Brownian motion $(B(u))_{0\leq u \leq 1}$ in $D[0,1]$ (see the Proof of Theorem 5 of \cite{RefAT}).
Since 
\[ \sup_{u\in[0,1]} \left| \sum_{j =1}^{\lfloor n^u \rfloor} \frac{\theta}{j} -  u \theta \log{n}  \right| = O(\theta), \]
the random process
\begin{equation}
 \left( \frac{ \sum_{j =1}^{\lfloor n^u \rfloor} C_j^\infty - u \theta \log{n}  }{ \sqrt{ \theta \log{n}  } } \right)_{0\leq u\leq1}    \label{FC2}
\end{equation}
converges weakly to $(B(u))_{0\leq u \leq 1}$ in $D[0,1]$ because of the assumption \eqref{FASS}.
From \eqref{FC1} and the weak convergence of \eqref{FC2}, Theorem 2.7 (iv) of \cite{RefV} yields the result.

(ii) First we argue $X^2_n(\cdot)$.
From \eqref{EB2} and \eqref{FASS1}, it follows that
\begin{equation*}
\int_0^1 \left|  \frac{\sum_{j=1}^{\lfloor n^u \rfloor}(C_j^n - C_j^\infty)}{\sqrt{ \sum_{j=1}^{\lfloor n^u \rfloor} \theta/j }} \right|^2 du \to^p 0.
\end{equation*}
It holds that
\begin{equation}
 \left( \frac{ \sum_{j =1}^{\lfloor n^u \rfloor} C_j^\infty - \sum_{j =1}^{\lfloor n^u \rfloor} \theta/j }{ \sqrt{ \sum_{j =1}^{\lfloor n^u \rfloor} \theta/j  } } \right)_{0<u<1} =^d  \left( \frac{ N^1(\sum_{j =1}^{\lfloor n^u \rfloor} \theta/j) - \sum_{j =1}^{\lfloor n^u \rfloor} \theta/j  }{ \sqrt{ \sum_{j =1}^{\lfloor n^u \rfloor} \theta/j } } \right)_{0<u<1} ,  \label{FC4}
\end{equation}
where $(N^1(t))_{t \geq 0}$ is a homogeneous Poisson process with unit intensity.
Since
\[ \sup_{u\in(0,1)} \frac{ | \sum_{j=1}^{ \lfloor n^u \rfloor} \theta/j  - u \theta \log{n} | }{ \theta \log{n}} \to 0\]
and the other hypotheses hold with $\lambda=1$, $s_n(u)= \sum_{j=1}^{\lfloor n^u \rfloor} \theta/j $ and $f(n)=\theta\log{n}$ (see Subsection 6.2 of \cite{RefT1}) Lemma \ref{PFCLT} in Appendix implies that \eqref{FC4} converges weakly to $(B(u)/\sqrt{u})_{0<u<1}$ in $L^2(0,1)$.
From what has been already proven, Theorem 2.7 (iv) of \cite{RefV} yields the result.

Next we argue $X^3_n(\cdot)$.
Since
\begin{eqnarray*}
&& \int_0^{\frac{\varepsilon}{\log{n}}} \frac{( N^1{(u \theta \log{n})} - u \theta \log{n} )^2}{u \theta \log{n}} du \to^p 0, \\
&& \int_{\frac{\varepsilon}{\log{n}}}^1 \frac{ \left\{ \sum_{j=1}^{\lfloor n^u \rfloor} (C_j^n - C_j^\infty)  \right\}^2}{u \theta \log{n}} du \to^p 0, \\
&& \int_{\frac{\varepsilon}{\log{n}}}^1 \frac{ \left( \sum_{j=1}^{\lfloor n^u \rfloor}C_j^\infty - N^1{(u \theta \log{n})}  \right)^2}{u \theta \log{n}} du \to^p 0 ,
\end{eqnarray*}
which follow from the almost same argument as the proof of Theorem 7.1 of \cite{RefT1} by the assumption \eqref{FASS1}, we have
\[ \int_0^1 \left( X^3_n(u) - \frac{N^1(u\theta \log{n}) - u\theta\log{n}}{\sqrt{u\theta\log{n}}} \right)^2 du \to^p 0. \]
From Lemma \ref{PFCLT} in Appendix with $\lambda=1$, $s_n(u) = u f(n)$ and $ f(n)= \theta \log{n}$,  it holds that the random process
\begin{equation*}
\left( \frac{N^1(u\theta \log{n}) - u\theta\log{n}}{\sqrt{u\theta\log{n}}} \right)_{0<u<1} 
\end{equation*}
converges weakly to $(B(u)/\sqrt{u})_{0<u<1}$ in $L^2(0,1)$.
Consequently, the desired result follows.
This completes the proof. \\
 
\begin{rem}
It follows from Proposition \ref{FCLT1} (i) that if \eqref{FASS} holds then \eqref{WCLT} holds.
But as it is stated in \eqref{TC}, the asymptotic normality of $K_n$ holds for far larger $\theta$.
\end{rem}
 
The following result, promised in Subsection \ref{fcs}, is an extension of Proposition \ref{LR0}.

\begin{prop}\label{LRC}
(i) In Case A, if \eqref{FASS} holds then the random process $X^4_n(\cdot)$ defined in \eqref{EB} converges weakly to a standard Brownian bridge $(B^\circ (u))_{0\leq u \leq 1}$ in $D[0,1]$. \\
(ii) In Case A, if \eqref{FASS1} holds then the random process $X^5_n(\cdot)$ defined in \eqref{EBT} converges weakly to $(B^\circ (u)/\sqrt{u(1-u)})_{0 < u < 1}$ in $L^2(0,1)$. 
\end{prop}

Proof.
(i) Since it holds that
\begin{eqnarray*}
X^4_n(u)  &=&  \frac{\theta \log{n}}{K_n} \left( \frac{\sum_{j=1}^{\lfloor n^u \rfloor} C_j^n - uK_n}{\sqrt{\theta \log{n}}} \right) \\
&=& \frac{\theta \log{n}}{K_n} \left\{  \frac{ (1-u) \sum_{j=1}^{\lfloor n^u \rfloor} C_j^n - u \sum_{j=\lfloor n^u \rfloor + 1}^{n} C_j^n}{\sqrt{\theta \log{n}}} \right\} 
\end{eqnarray*}
for any $u \in [0,1]$, it is sufficient to show
\begin{equation}
\frac{K_n}{\theta \log{n}} \to^p  1 \label{lrc2}
\end{equation}
and
\begin{equation}
\left(  \frac{ (1-u) \sum_{j=1}^{\lfloor n^u \rfloor} C_j^n - u \sum_{j=\lfloor n^u \rfloor + 1}^{n} C_j^n}{\sqrt{\theta \log{n}}}  \right)_{0\leq u \leq 1} \Rightarrow (B^\circ(u))_{0\leq u \leq 1}  \label{lrc1} 
\end{equation}
in $D[0,1]$.
Firstly, \eqref{lrc2} holds because the assumption \eqref{FASS} yields that $\log{\theta}/\log{n} \to 0$ and because it follows from \eqref{FLLN} that
\[ \frac{K_n}{\theta\log{n} } = \frac{K_n}{\theta\log{(n/\theta)} \left(1 + \frac{\log{\theta}}{\log{(n/\theta)}} \right) } \to^p 1. \]
Next, we show \eqref{lrc1}.
Since it follows from Proposition \ref{PPAU} and from the assumption \eqref{FASS} (see \eqref{FC1}) that
\[\sup_{u\in[0,1]} \left|  \frac{ \sum_{j=1}^{\lfloor n^u \rfloor} C_j^n - \sum_{j=1}^{\lfloor n^u \rfloor} C_j^\infty }{\sqrt{\theta\log{n}}} \right| \leq \frac{\sum_{j=1}^n |C_j^n - C_j^\infty|}{\sqrt{\theta \log{n}}} \to^p 0 \]
and that
\[ \sup_{u\in[0,1]} \left|  \frac{ \sum_{j=\lfloor n^u \rfloor + 1}^{n} C_j^n - \sum_{j=\lfloor n^u \rfloor + 1}^{n} C_j^\infty }{\sqrt{\theta\log{n}}} \right|\leq \frac{\sum_{j=1}^n |C_j^n - C_j^\infty|}{\sqrt{\theta \log{n}}} \to^p 0, \]
the triangle inequality yields that
\[ \sup_{u\in[0,1]} \left|  \frac{ (1-u) \sum_{j=1}^{\lfloor n^u \rfloor} C_j^n - u \sum_{j=\lfloor n^u \rfloor + 1}^{n} C_j^n}{\sqrt{\theta \log{n}}}  -  P_4^\circ(u)  \right| \to^p 0  \]
where 
\[ \left( P_4^\circ(u) \right)_{0 \leq u \leq 1} = \left(  \frac{ \sum_{j=1}^{\lfloor n^u \rfloor} C_j^\infty - u \sum_{j=1}^{n} C_j^ \infty }{\sqrt{\theta \log{n}}} \right)_{0\leq u \leq 1}. \]
By using the functional central limit theorem for Poisson processes in $D[0,1]$, $P_4^\circ(\cdot)$ converges weakly to $(B^\circ (u))_{0\leq u \leq 1}$ in $D[0,1]$ (see \cite{RefAT}).

(ii) By the same reason as (i), it is sufficient to show \eqref{lrc2} and
\begin{eqnarray}
&& \left(  \frac{ (1-u) \sum_{j=1}^{\lfloor n^u \rfloor} C_j^n - u \sum_{j=\lfloor n^u \rfloor + 1}^{n} C_j^n}{\sqrt{u(1-u) \theta \log{n}}} 1\left\{ \frac{\varepsilon}{\log{n}} < u < 1-\frac{\varepsilon}{\log{n}} \right\} \right)_{0 < u < 1} \nonumber \\
&\Rightarrow& \left(\frac{B^\circ(u)}{\sqrt{u(1-u)}} \right)_{0 < u < 1}  \label{lrc3} 
\end{eqnarray}
in $L^2(0,1)$.
Firstly, it holds that
\begin{eqnarray}
&& \int_{\frac{\varepsilon}{\log{n}}}^{1- \frac{\varepsilon}{\log{n}} } \frac{(1-u)^2 \left\{ \sum_{j=1}^{\lfloor n^u \rfloor} (C_j^n - C_j^\infty)  \right\}^2 }{u(1-u) \theta \log{n}} du \nonumber \\
&\leq& \int_{\frac{\varepsilon}{\log{n}}}^{1 } \frac{ \left( \sum_{j=1}^{\lfloor n^u \rfloor} |C_j^n - C_j^\infty |  \right)^2}{u \theta \log{n}} du \nonumber  \\
&\leq& \int_{\frac{\varepsilon}{\log{n}}}^{1 } \frac{ \left( \sum_{j=1}^{n} |C_j^n - C_j^\infty |  \right)^2}{u \theta \log{n}} du \nonumber  \\
&\leq& \frac{\log\log{n} - \log{\varepsilon} }{\theta\log{n}} \left( \sum_{j=1}^{n} |C_j^n - C_j^\infty |  \right)^2 \label{prp631},
\end{eqnarray}
and that
\begin{eqnarray}
&& \int_{\frac{\varepsilon}{\log{n}}}^{1- \frac{\varepsilon}{\log{n}} } \frac{ u^2 \left\{ \sum_{j=\lfloor n^u \rfloor+1}^{n} (C_j^n - C_j^\infty)  \right\}^2}{u(1-u) \theta \log{n}} du \nonumber  \\
&\leq&  \int_{0}^{1- \frac{\varepsilon}{\log{n}} } \frac{ \left( \sum_{j=\lfloor n^u \rfloor+1}^{n} |C_j^n - C_j^\infty|  \right)^2}{(1-u) \theta \log{n}} du \nonumber  \\
&\leq&  \frac{\log\log{n} - \log{\varepsilon} }{\theta\log{n}} \left( \sum_{j=1}^{n} |C_j^n - C_j^\infty |  \right)^2 \label{prp632}.
\end{eqnarray}
These right-hand sides of \eqref{prp631} and \eqref{prp632} converge to 0 in probability because the expectations of their square root converge to 0 by the assumption \eqref{FASS1}.
Secondly, by letting $(N^1(t))_{t \geq 0}$ be a homogeneous Poisson process with unit intensity, it follows from
\begin{eqnarray*}
{\sf E}[( N^1(u \theta \log{n}) - u N^1(\theta \log{n}) )^2] &=& {\rm var}(N^1(u \theta \log{n}) - u N^1(\theta \log{n})) \\
&=& (1-u)^2 u \theta \log{n} + u^2 (1-u) \theta \log{n} \\
&=& u(1-u)\theta\log{n}
\end{eqnarray*}
that
\[
\int_0^{\frac{\varepsilon}{\log{n}}} \frac{( N^1(u \theta \log{n}) - u N^1(\theta \log{n}) )^2}{u(1-u) \theta \log{n}} du \to^p 0, \]
and that
\[
 \int_{1-\frac{\varepsilon}{\log{n}}}^1 \frac{( N^1(u \theta \log{n}) - u N^1(\theta \log{n}) )^2}{u(1-u) \theta \log{n}} du \to^p 0.
\]
Thirdly, it holds that
\begin{eqnarray*}
&& \int_{\frac{\varepsilon}{\log{n}}}^{1- \frac{\varepsilon}{\log{n}}} \frac{ \left\{ \sum_{j=1}^{\lfloor n^u \rfloor}C_j^\infty - u \sum_{j=1}^{n }C_j^\infty - (N^1(u \theta \log{n}) - u N^1(\theta \log{n})  ) \right\}^2}{u(1-u) \theta \log{n}} du\\
&=&  \int_{\frac{\varepsilon}{\log{n}}}^{1- \frac{\varepsilon}{\log{n}}} 
\frac{1}{u(1-u) \theta \log{n}} 
 \left[ (1-u) \left( \sum_{j=1}^{\lfloor n^u \rfloor} C_j^\infty - N^1(u \theta \log{n}) \right) \right. \\ && \left. - u \left\{ \sum_{j=\lfloor n^u \rfloor +1}^{n }C_j^\infty - \left( N^1(\theta \log{n})- N^1(u \theta \log{n}) \right)  \right\} \right]^2 du \\
&\leq& \int_{\frac{\varepsilon}{\log{n}}}^{1- \frac{\varepsilon}{\log{n}}} \frac{2 (1-u)^2 \left( \sum_{j=1}^{\lfloor n^u \rfloor}C_j^\infty  - N^1(u \theta \log{n} ) \right)^2}{u(1-u) \theta \log{n}} du \\
&& \quad  + \int_{\frac{\varepsilon}{\log{n}}}^{1- \frac{\varepsilon}{\log{n}}} \frac{2 u^2 \left\{ \sum_{j=\lfloor n^u \rfloor +1}^{n }C_j^\infty -  (N^1(\theta \log{n}) - N^1(u \theta \log{n})  ) \right\}^2}{u(1-u) \theta \log{n}} du \\
&\leq& \int_{\frac{\varepsilon}{\log{n}}}^{1} \frac{2 \left( \sum_{j=1}^{\lfloor n^u \rfloor}C_j^\infty  - N^1(u \theta \log{n} ) \right)^2}{u\theta \log{n}} du \\
&& \quad  + \int_{0}^{1- \frac{\varepsilon}{\log{n}}} \frac{ 2 \left\{ \sum_{j=\lfloor n^u \rfloor +1}^{n }C_j^\infty -  (N^1(\theta \log{n}) - N^1(u \theta \log{n})  ) \right\}^2}{(1-u) \theta \log{n}} du.
\end{eqnarray*}
The distributions of the first term and second term in the right-hand side are equal to
\[ \int_{\frac{\varepsilon}{\log{n}}}^{1} \frac{2 \left(  N^1(\theta( \sum_{j=1}^{\lfloor n^u \rfloor}  1/j - u \log{n}) ) \right)^2}{u\theta \log{n}} du \]
and
\[ \int_{0}^{1- \frac{\varepsilon}{\log{n}}} \frac{2 \left( N^1(\theta( \sum_{j=1}^{n} 1/j - \log{n}) ) - N^1(\theta( \sum_{j=1}^{\lfloor n^u \rfloor}  1/j - u \log{n}) ) ) \right)^2}{(1-u) \theta \log{n}} du, \]
respectively.
Both of them converge to 0 in probability because their expectations tend to 0 from the assumption \eqref{FASS1}.
Thus, the triangle inequality yields that
\[ \int_0^1 \left|  X_n^5(u) -  P_5^\circ(u)  \right|^2 du \to^p 0,  \]
where 
\[ \left( P_5^\circ(u) \right)_{0 < u < 1} 
= \left(  \frac{ N^1(\theta u \log{n}) - u N^1(\theta \log{n})  }{\sqrt{u(1-u) \theta \log{n}}}  \right)_{0 < u < 1}.\]
Since
\begin{eqnarray*}
\left( P_5^\circ(u) \right)_{0 < u < 1} 
&=&  \left(  \int_0^{\theta \log{n}} \frac{ 1\{ t\leq u \theta \log{n} \} -u  }{\sqrt{u(1-u) \theta \log{n}}} dN^1(t) \right)_{0 < u < 1} \\
&=&  \left(  \int_0^{\theta \log{n}} \frac{ 1\{ t\leq u \theta \log{n} \} -u  }{\sqrt{u(1-u) \theta \log{n}}} (dN^1(t)  - dt) \right)_{0 < u < 1} ,
\end{eqnarray*}
Theorem 4 of \cite{RefT16} yields that $\left( P_5^\circ(u) \right)_{0 < u < 1}  \Rightarrow \left( B^\circ(u)/\sqrt{u} \right)_{0 < u < 1} $ by setting $H_s = 1$ with $d=1$, $\lambda_s=1$ and $T=\theta\log{n}$.
Consequently, we have \eqref{lrc3}.
This completes the proof.

%%%%%%%%%%%%%%%%%%%%%%%%%%%%%%%%%%%%%%%%%%%%%%%%%%%%%

\appendix
\section{Appendix. Auxiliary lemmas}

\begin{lem}\label{NBL}
In Case A, if $\theta^2/n \to0$ then
\begin{equation}
\Gamma(\theta)  \frac{(\theta)_n }{n!} =  n^{\theta-1} \left\{ 1+ \frac{\theta(\theta -1 )}{2n} \right\} +O \left( n^{\theta-1} \frac{\theta^4}{n^2} \right).  \label{la1s}
\end{equation}
\end{lem}

Proof.
The left-hand side of \eqref{la1s} equals 
\[\frac{\Gamma(n+ \theta)}{n\Gamma(n)}. \]
By using the asymptotic series expansion $\Gamma(x) = \sqrt{2\pi} e^{-x}x^{x-1/2} \{ 1+1/(12x) \} + O\left( e^{-x}x^{x-5/2}  \right)$ as $x\to \infty$, it holds that
\begin{eqnarray*}
\Gamma (n+\theta)
&=& \sqrt{2\pi} e^{-(n+ \theta)} (n+ \theta)^{n+\theta-1/2} \left\{  1+ \frac{1}{12(n+ \theta)} \right\} +O \left( e^{-(n+ \theta)} (n+ \theta)^{n+\theta-5/2}  \right)   \\
&=& \sqrt{2\pi} e^{-(n+ \theta)} n^{n+\theta-1/2}  \left( 1+ \frac{\theta}{n} \right)^{n+\theta-1/2} \left(  1+ \frac{1}{12n} + O\left( \frac{\theta}{n^2} \right) \right) 
\end{eqnarray*}
and that
\begin{eqnarray*}
n \Gamma(n)
&=& n \sqrt{2\pi} e^{-n} n^{n-1/2} \left(  1+ \frac{1}{12n} \right) +O \left( e^{-n} n^{n-5/2}  \right) \\
&=& \sqrt{2\pi} e^{-n} n^{n+1/2} \left(  1+ \frac{1}{12n} \right) \left( 1+ O \left( \frac{1}{n^2} \right) \right) .
\end{eqnarray*}
Hence, the left-hand side of \eqref{la1s} is
\begin{eqnarray}
&& n^{\theta-1} e^{-\theta} \left( 1 + \frac{\theta}{n} \right)^n \left( 1 + \frac{\theta}{n} \right)^{\theta- 1/2} \left( 1+ O\left( \frac{\theta}{n^2} \right) \right) \left( 1+O\left( \frac{1}{n^2} \right) \right) \nonumber \\
&=& n^{\theta-1} e^{-\theta} \left( 1 + \frac{\theta}{n} \right)^n \left( 1 + \frac{\theta}{n} \right)^{\theta- 1/2} \left( 1+ O\left( \frac{\theta}{n^2} \right) \right)  \label{tfG}
\end{eqnarray}
and, from the asymptotic expansion $\left(1+ 1/x \right)^x = e \left( 1- 1/(2x) +O\left( x^{-2} \right) \right)$ as $x\to\infty$ it follows that 
\[\left( 1 + \frac{\theta}{n} \right)^n = e^\theta \left( 1- \frac{\theta}{2n} +O\left( \frac{\theta^2}{n^2} \right) \right)^\theta\]
 and hence \eqref{tfG} is
\begin{eqnarray*}
&&  n^{\theta-1} \left( 1 - \frac{\theta}{2n} +O\left( \frac{\theta^2}{n^2} \right) \right)^\theta \left( 1 + \frac{\theta}{n}  \right)^\theta  \left( 1 + \frac{\theta}{n} \right)^{- 1/2}  \left( 1+ O\left( \frac{\theta}{n^2} \right) \right)  \\
&=& n^{\theta-1} \left( 1 + \frac{\theta}{2n} +O\left( \frac{\theta^2}{n^2} \right) \right)^\theta   \left( 1 + \frac{\theta}{n} \right)^{- 1/2}  \left( 1+ O\left( \frac{\theta}{n^2} \right) \right)  \\
&=& n^{\theta-1}   \left( 1 + \frac{\theta^2}{2n} +O\left( \frac{\theta^4}{n^2} \right) \right) \left( 1 - \frac{\theta}{2n} + O\left( \frac{\theta^2}{n^2} \right) \right) \left( 1+ O\left( \frac{\theta}{n^2} \right) \right) \\
&=& n^{\theta-1}   \left\{ 1 + \frac{\theta(\theta-1)}{2n}   +O\left( \frac{\theta^4}{n^2} \right) \right\}  \left( 1+ O\left( \frac{\theta}{n^2} \right) \right) \\
&=&  n^{\theta-1}   \left\{ 1 + \frac{\theta(\theta-1)}{2n}   +O\left( \frac{\theta^4}{n^2} \right) \right\} .
\end{eqnarray*}
This completes the proof.

\begin{lem}\label{SI}
Let $b$ be an positive integer.
For $a>1$, it holds that
\[ \sum_{j=1}^b \frac{a^j}{j} \leq \log{b} + a^b.  \]
\end{lem}

Proof.
It holds that
\[  \sum_{j=1}^b \frac{a^j}{j}  =  \sum_{j=1}^b \int_0^a t^{j-1} dt =  \sum_{j=1}^b \left( \frac{1}{j} +  \int_1^a t^{j-1} dt \right) \leq 1+\log{b} + \int_1^a \frac{t^b-1}{t-1} dt. \]
As for the last term, it holds that
\begin{eqnarray*}
\int_1^a \frac{t^b-1}{t-1} dt 
&=& \int_1^a \frac{ \{ 1+(t-1) \}^b  - 1}{t-1} dt \\
&=& \int_1^a \sum_{j=1}^b \binom{b}{j} (t-1)^{j-1} dt 
=  \int_1^a \sum_{j=1}^b \binom{b-1}{j-1} \frac{b}{j} (t-1)^{j-1}  dt \\
&\leq&  \int_1^a b \sum_{j=1}^b \binom{b-1}{j-1} (t-1)^{j-1} dt 
=   \int_1^a b t^{b-1} dt \\
&=& a^b -1.
\end{eqnarray*}
This completes the proof.

\begin{lem}\label{SIL}
For any $a>0$ and for any positive integer $b$, $(x-a)_b/(x)_b$ is increasing with respect to $x>a$.
\end{lem}

Proof.
The proof is by induction on $b$.
When $b=1$, $ (x-a)/x = 1 - a/ x$ is increasing.
Let $x_1$ and $x_2$ satisfy $a < x_1 <x_2$.
If the conclusion of the lemma is true for $b$, then the conclusion is also true for $b+1$ because
\[ \frac{(x_2-a)_{b+1}}{(x_2)_{b+1} } =  \frac{x_2+b-a}{x_2+b} \frac{(x_2-a)_b}{(x_2)_b} > \frac{x_2+b-a}{x_2+b} \frac{(x_1-a)_b}{(x_1)_b}  > \frac{x_1+b-a}{x_1+b} \frac{(x_1-a)_{b}}{(x_1)_{b} }  =  \frac{(x_1-a)_{b+1}}{(x_1)_{b+1} } . \]
This completes the proof.

\begin{lem}\label{PFCLT}
Let $( N_t )_{t\geq 0}$ be a homogeneous Poisson process with intensity $\lambda>0$ satisfying $N_0=0$.
Define the non-decreasing function $(n,u) \mapsto s_n(u)$ with respect to $0 \leq u \leq 1$ and with respect to $n=1,2,\ldots$ which satisfies
$ \inf_{u \in (\tau , 1)} s_n(u) >0 $
for all $0 < \tau <1$,
\begin{equation}
\lim_{n\to\infty} \left( \frac{ \sup_{u \in (0,1)} \left| s_n(u) - u f(n) \right|}{f(n)}  \right) = 0  \label{Pois1}
\end{equation}
with $n \mapsto f(n)$ an increasing function of $n$ satisfying $\lim_{n\to\infty} f(n) = \infty $, and
\[
\lim_{n\to\infty} \left\{ \int_0^1 \frac{du}{(s_n(u))^\delta} \right\} = 0 
\]
for some $\delta>0$.
Then, the random process
\[ \left( \frac{N_{s_n(u)} - \lambda s_n(u)}{\sqrt{\lambda s_n(u)}} \right)_{0<u<1} \]
converges weakly to a Gaussian process $(B(u)/\sqrt{u})_{0<u<1}$ in $L^2(0,1)$ as $n \to \infty$.
\end{lem}

\begin{rem}
Lemma \ref{PFCLT} is a slight generalization of Lemma 2.1 of \cite{RefT1}.
The only difference is condition \eqref{Pois1}, where corresponding condition (2.1) of \cite{RefT1} is the case that $f(n)=K \log{n}$ with a constant $K$.
To show Lemma \ref{PFCLT}, the equation
\[ \lim_{n\to\infty} \left( \frac{s_n(u) \wedge s_n(v) }{\sqrt{s_n(u)s_n(v)}} \right) = \lim_{n\to\infty} \left( \frac{K ((u \log n ) \wedge (v \log n))}{\sqrt{K^2 \log n^u \log n^v}} \right) =\frac{u\wedge v}{\sqrt{uv}} \]
in the proof of Lemma 2.1 of \cite{RefT1} should be replaced by
\[ \lim_{n\to\infty} \left( \frac{s_n(u) \wedge s_n(v) }{\sqrt{s_n(u)s_n(v)}} \right) = \lim_{n\to\infty} \left\{ \frac{ (u f(n) ) \wedge (v f(n)) }{\sqrt{ u f(n) v f(n)}} \right\} =\frac{u\wedge v}{\sqrt{uv}}, \]
and the other part has no need to change.
\end{rem}

%%%%%%%%%%%%%%%%%%%%%%%%%%%%%%%%%%%%%%%%%%%%%%%%%%%%%%%%%%%
\section*{Acknowledgements}
The author would like to his heartfelt gratitude to Professor Shuhei Mano for a lot of constructive comments.
%%%%%%%%%%%%%%%%%%%%%%%%%%%%%%%%%%%%%%%%%%%%%%%%%%%%%%%%%%%

\end{document}